\documentclass[aps,twocolumn,preprintnumbers]{revtex4-2}

\usepackage{graphicx}  
\usepackage[caption=false]{subfig}
\usepackage{multirow}
\usepackage{fancyhdr}
\usepackage{longtable}
\usepackage{parskip}
\usepackage[T1]{fontenc}
\usepackage{dcolumn}   

\usepackage{bm}        
\usepackage{amsfonts}  
\usepackage{amsmath}   
\usepackage{amssymb}   
\usepackage{hyperref}

\hypersetup{
    colorlinks=true,
    linkcolor=cyan,
    citecolor=cyan,
    urlcolor=cyan}

\usepackage{float}

\newcommand{\pwisein}{\left\{ \begin{array}{ll}}
\newcommand{\pwiseout}{\end{array}\right.}

\newcounter{defcounter}
\usepackage{amsthm}

\begin{document}
\title[]{Computing leaky Lamb waves for waveguides between elastic half-spaces using spectral collocation}

\author{Evripides Georgiades}
\email{evripides.georgiades18@imperial.ac.uk}
\author{Michael J.S. Lowe}

\author{Richard V. Craster}
\altaffiliation{Also at: Department of Mathematics, Imperial College London, London, SW7 1AY, United Kingdom}

\affiliation{Department of Mechanical Engineering, Imperial College London, London, SW7 1AY, United Kingdom}

\date{\today} 

\begin{abstract}

In non-destructive evaluation guided wave inspections, the elastic structure to be inspected is often embedded within other elastic media and the ensuing leaky waves are complex and non-trivial to compute; we consider the canonical example of an elastic waveguide surrounded by other elastic materials that demonstrates the fundamental issues with calculating the leaky waves in such systems. Due to the complex wavenumber solutions required to represent them, leaky waves pose significant challenges to existing numerical methods, with methods that spatially discretise the field to retrieve them suffering from the exponential growth of their amplitude far into the surrounding media. We present a spectral collocation method yielding an accurate and efficient identification of these modes, leaking into elastic half-spaces. We discretise the elastic domains and, depending on the exterior bulk wavespeeds, select appropriate mappings of the discretised domain to complex paths, in which the numerical solution decays and the physics of the problem are preserved. By iterating through all possible radiation cases, the full set of dispersion and attenuation curves are successfully retrieved and validated, where possible, against the commercially available software DISPERSE. As an independent validation, dispersion curves are obtained from finite element simulations of time-dependent waves using Fourier analysis. 

\end{abstract}

\maketitle

\section{Introduction}

 Computing leaky guided waves, and in particular accurately identifying their attenuation and dispersion, is crucial for developing effective practical ultrasonic inspection techniques that use guided waves for non-destructive evaluation (NDE), e.g. in the assessment of structural integrity. In the context of interest here, the structure under inspection comprises an elastic waveguide in contact with dissimilar elastic media; energy is no longer necessarily confined within the boundaries of a waveguide but rather leaks into the adjacent media and causes the guided wave to attenuate. This behaviour is common in a wide range of settings, some of which extend beyond structural integrity, including reinforced concrete structures \citep{Shen2014DispersionStudy}, seismology \citep{Lellouch2021PropertiesFiber}, adhesive joints and fiber composites \citep{Chimenti1985LeakyLaminates,Nagy1994LeakyMaterials}, buried pipes \citep{Long2003AttenuationPipes,Long2004MeasurementWaveguide,Leinov2016UltrasonicPipes,Leinov2015InvestigationSand,Castaings2008FiniteMedia} and biological structures such as bones surrounded by tissue \citep{Potsika2014ApplicationBones,Guha2021IdentificationModelling,Lee2004FeasibilityTibia}. Consequently, understanding the dispersion and attenuation characteristics of leaky systems, and in particular leaky Lamb waves \citep{Lamb1916OnPlate,Kiefer2019CalculatingInteraction, Pavlakovic1998LeakyNDT, Georgiades2022LeakyMethods}, for waveguides in contact with spaces of different elastic materials is of recurrent and broad interest.  

The modal solutions of leaky Lamb waves often display the amplitude of the leaked energy exponentially growing with distance into the adjacent elastic media \citep{Treyssede2014FiniteWaveguides,Chimenti1985LeakyLaminates,Kiefer2019CalculatingInteraction}. This is a consequence of energy conservation and this phenomenon is discussed in detail\citep{Georgiades2022LeakyMethods} for the case of acoustic energy radiation from an elastic waveguide immersed in fluid. The exponential growth in amplitude and the complex wavenumber required to describe leaky wave modes, both pose a significant challenge for conventional numerical solvers used to model guided elastic wave propagation. For instance, Partial Wave Root Finding methods (PWRF) are well suited for solving wave propagation problems of real wavenumber solutions, i.e. non-radiating problems, however they struggle to converge to the attenuated waves characterised by a complex wavenumber \citep{Rokhlin1989OnLayer}. Other methods, such as Semi Analytical Finite Element models, require elaborate techniques to solve the nonlinear eigenvalue problems \citep{Hayashi2014CalculationMethod,Mazzotti2014UltrasonicValidation,Mazzotti2013ACross-section} that emerge or require a carefully defined Perfectly Matched Layer (PML) to accurately model wave propagation in an infinite medium \citep{Treyssede2014FiniteWaveguides,Zuo2017NumericalPackage}.
 
An alternative to traditional methods is the use of Spectral Collocation Methods\citep{Boyd2000ChebyshevMethods,Trefethen2008SpectralMatlab}  (SCMs) for the discretisation of the solution domains, while SCM has had considerable success for magnetoelastic \citep{Ryskamp2023ComputationMaterials} or elastic guided waves \citep{Adamou2004SpectralMedia} in layered media, anisotropic and viscoelastic media \citep{Quintanilla2016FullMedia,Quintanilla2015GuidedMedia,Quintanilla2015ModelingMethod,Quintanilla2017TheWaveguides}. Recently, for a waveguide immersed with an inviscid fluid, Kiefer \textit{et. al } \citep{Kiefer2019CalculatingInteraction} employed an SCM,  and exact interface conditions, to account for the infinite fluid domain in contact with the waveguide; the non-linear eigenvalue problem obtained from this discretisation was subsequently linearised using a change of variable. In a similar manner, Tang \textit{et. al} \citep{Tang2022StudyFluids} created a non-linear eigenvalue problem for a waveguide in contact with a different fluid on either side, which they then linearised into a quintic polynomial eigenvalue problem by moving to a higher dimensional state space: It would be a major challenge to attempt to implement this in elasticity for a different medium in each side of the waveguide. 

Using an alternative approach \cite{Georgiades2022LeakyMethods}, we identified the leaky Lamb waves emitted from an infinite waveguide fully immersed in inviscid fluids also using SCMs and overcame the challenges of an exponentially growing solution by mapping the discretised fluid domains to complex paths, in which the numerical solution decayed, while maintaining the physics of wave propagation in the physical space. Here, we extend to the waveguide not in contact with a fluid but rather with elastic half-spaces. As elastic materials support both longitudinal and shear waves, it is necessary to implement a suitable discretisation approach to accommodate both these wave types; this is a step-change in complexity from fluid loading cases, as either one or both longitudinal/shear waves can exponentially grow for leaky waves. The key contribution of this work will be the development of a methodology, based on the solution method of \citep{Georgiades2022LeakyMethods}, that accounts for the presence of two distinct partial waves within each exterior elastic domain and solves the challenging problem of leakage into solids - a problem of practical relevance to the NDE community. We will demonstrate a SCM, in which the discretisation of the exterior elastic domains is determined by the nature of the solution to be found. This method delivers complete spectra of the dispersion curves and their mode shapes, with high accuracy, and as ``black box" solutions: there is no need for priori knowledge of the nature of the solutions nor for the tuning of parameters to approach the solution. 

To validate our approach quantitatively, we will firstly compare our results with those obtained using the commercially available PWRF software DISPERSE \citep{Pavlakovic1997Disperse:Curves}. We will then strengthen our validation by comparing our findings with Finite Element (FE) modelling results. As FE simulations are agnostic to the modal concepts introduced by our method and only result in a weighted superposition of modes\citep{BatheK.J.1996FiniteProcedures}, the comparison with our method will further enhance the credibility of our findings and will instill confidence in the robustness of our approach, particularly in regions where DISPERSE may not yield viable solutions.

This paper is organised as follows - Section \ref{Theory} presents a summary of the underlying theory of leaky Lamb waves radiating from a waveguide in contact with two different isotropic elastic half-spaces. Following that, Section \ref{Method} presents our discretisation method and FE model, while a discussion and the validation of our results can be found in Section \ref{Results and Discussion}.

\section{Theory}\label{Theory}
\begin{figure}
    \centering
    \includegraphics[width=\columnwidth]{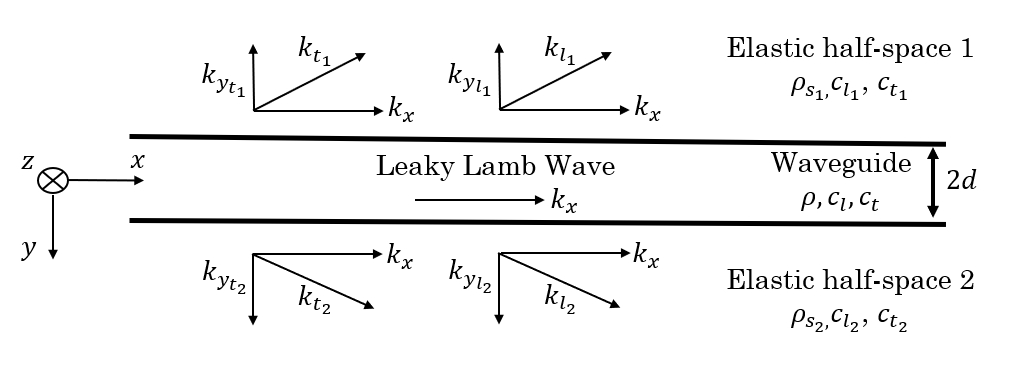}
    \caption{Schematic of an elastic waveguide, between two elastic half-spaces, of thickness $2d$ and infinite extent in the $x$ and $z$ directions. The densities of the elastic media of the setup are $\rho$, $\rho_{1}$ and $\rho_{2}$, their longitudinal wavespeeds are $c_l$, $c_{l_1}$ and $c_{l_2}$ and their transverse wavespeeds are $c_t$, $c_{t_1}$ and $c_{t_2}$ respectively. }
    \label{Schematic}
\end{figure}

We consider an isotropic, linear-elastic waveguide, of finite thickness, $2d$, and infinite extent in both the $x$ and $z$ directions, that is between two half-spaces of different isotropic and linear-elastic media. Both the elastic waveguide and half-spaces are assumed to be homogeneous with densities $\rho$, $\rho_{1}$ and $\rho_{2}$, longitudinal wavespeeds $c_l$, $c_{l_1}$ and $c_{l_2}$ and transverse wavespeeds $c_t$, $c_{t_1}$ and $c_{t_2}$ respectively. A pictorial representation of the cross section of the setup together with the axis configuration is shown in Fig. \ref{Schematic}. Consistent with the physical interpretation of the classical definition of Auld\citep{Auld1973AcousticSolidsb} for a free plate, we refer to those modes of non-zero energy flux as propagating modes; our investigation will focus on modes of positive phase velocities.

Energy radiation into the adjacent half-spaces is expressed mathematically in terms of a common complex wavenumber of propagation, denoted here by $k_x$, whose real part reveals information about the velocity of propagation and whose imaginary part corresponds to the attenuation of the propagated waves. In each of the solid half-spaces, the radiated waves are regarded as the superposition of two partial waves, one longitudinal and one shear, whose complex wavenumbers along the $y$-direction are $k_{y_{l_j}}$ and $k_{y_{t_j}}$ respectively, for $j=1,2$ for each of the two elastic half-spaces. Under this regime, for a given angular frequency, $\omega$, the following holds true \citep{Nayfeh1997ExcessLoading}:

\begin{equation} \label{wavenumber decompositiono}
    k_{y_{l_j}}^2=k_{l_j}^2-k_x^2, \text{ } k_{y_{t_j}}^2=k_{t_j}^2-k_x^2, \text{ } k_{y_{l_j}}=\frac{\omega}{c_{l_j}} \text{ and } k_{y_{t_j}}=\frac{\omega}{c_{t_j}}.
\end{equation}

\noindent These propagating and radiating waves, with displacements in only the $x$ and $y$ directions, are the leaky Lamb waves \citep{Georgiades2022LeakyMethods} we seek.  

To represent displacement fields, we denote by $\boldsymbol{\bar{u}}$ the vector displacement in the waveguide
\begin{equation}\label{displacement vector waveguide} 
        \boldsymbol {\bar{u}}=\begin{pmatrix}
\bar{u}_x \\
\bar{u}_y \\
0 \end{pmatrix},
\end{equation}
whereby $\bar{u}_x$ and $\bar{u}_y$ are the displacements in the propagation and transverse directions. Similarly, for the $x-y$ plane displacement fields in the two elastic half-spaces, we employ a Helmholtz decomposition \citep{Achenbach1973WaveSolids,Rose2014UltrasonicMedia}, expressing the vector displacements in the two media, $\boldsymbol{\bar{u}}_{1}$,  $\boldsymbol{\bar{u}}_{2}$ 
as 
\begin{equation} \label{displacement vectors half spaces} 
  \boldsymbol {\bar{u}}_{j}=\begin{pmatrix}
\bar{u}_{x_{j}} \\
\bar{u}_{y_{j}} \\
0 \end{pmatrix}=\nabla \bar{\phi}_{j} + \nabla \times \begin{pmatrix} 0 \\
0 \\
\bar{\Psi}_{j}\end{pmatrix}.
\end{equation}

\noindent In this decomposition, $\bar{\phi}_{j}$ 
are the longitudinal wave potentials in the spaces while $\bar{\Psi}_{j}$ designates the $z-$component of their transverse vector wave potentials, which we henceforth refer to as their transverse wave potentials. 

Each of $\bar{u}_x, \bar{u}_y, \bar{\phi}_{j}$ and $\bar{\Psi}_{j}$ are functions of $x,y$ and $t$, while, when plane harmonic waves are assumed to be propagating through the structure, the common exponential term $\exp\left[i(k_xx-\omega t)\right]$ fully describes their harmonic nature. Thus, they decompose into
\begin{equation}
\label{u_x} \bar{\phi}_{j}(x,y,t)={\phi}_{j}(y)\exp\left[i(k_xx-\omega t)\right],
\end{equation}
with ${\phi}_{j}(y)$ representing the $y-$dependent components of the longitudinal wave potentials in the half-spaces and similarly for $\bar{u}_x (x,y,t), \bar{u}_y (x,y,t)$ and $\bar{\Psi}_{j} (x,y,t)$.

With this notation, the governing equations of motion, resulting from momentum equations \citep{Georgiades2022LeakyMethods,Nayfeh1997ExcessLoading}, are given by 

\begin{align}
\label{equation of motion waveguide x}
    \rho \omega^2 u_x-(\lambda+2\mu)k_x^2u_x+(\lambda+\mu) ik_x\frac{\partial u_y}{\partial y}+\mu\frac{\partial^2 u_x}{\partial y^2}&=0, \\
\label{equation of motion waveguide y}
    \rho \omega^2 u_y-\mu k_x^2u_y+(\lambda+\mu) ik_x\frac{\partial u_x}{\partial y}+(\lambda+2\mu)\frac{\partial^2 u_y}{\partial y^2}&=0, \\
\label{equation of motion half-space longitudinal}
\rho_j\omega^2\phi_{j}-k_x^2(\lambda_{j}+2\mu_{j})\phi_{j}+(\lambda_{j}+2\mu_{j})\frac{\partial^2 \phi_{j}}{\partial y^2}&=0, \\
\label{equation of motion half-space shear}
    \rho_{j}\omega^2\Psi_{j}-k_x^2\mu_{j}\Psi_{j}+\mu_{j}\frac{\partial^2 \Psi_{j}}{\partial y^2} &=0, 
\end{align}

\noindent with $\lambda$ and $\mu$, and $\lambda_{j}$ and $\mu_{j}$, the elastic Lam\'e parameters of the waveguide and the two elastic half-spaces respectively. 

The physics of the problem require the continuity of tractions and displacements at the two interfaces. That is to say, that at the two interfaces $y_1=-d$ and $y_2=d$, $\sigma_{y}-\sigma_{y_{j}}=0, \sigma_{xy}-\sigma_{xy_{j}}=0, u_x-u_{x_{j}}=0$ and $u_y-u_{y_{j}}=0$, with $\sigma_{y}$ and $\sigma_{y_j}$ the normal and $\sigma_{xy}$ and $\sigma_{xy_j}$ the shear stresses in the waveguide and half-spaces respectively \citep{Nayfeh1997ExcessLoading,Achenbach1973WaveSolids}.

Explicitly, the continuity of normal and shear stresses at the two interfaces is given by Eqs. \eqref{continuity of normal stress} and \eqref{continuity of shear stress} in terms of wave potentials in the elastic half-spaces

\begin{align}
\label{continuity of normal stress} 
&\left. \left[(\lambda+2\mu)\frac{\partial u_y}{\partial y}+\lambda ik_x u_x+\lambda_j k_x^2 \phi_j \right. \right. \nonumber \\ 
& \left. \left. \quad -(\lambda_j+2\mu_j)\frac{\partial^2 \phi_j}{\partial y^2}-2\mu_j ik_x \frac{\partial \Psi_j}{\partial y}\right]\right|_{y=y_j}=0, \\
\label{continuity of shear stress} 
&\left. \left[\mu \left( ik_x u_y+\frac{\partial u_x}{\partial y}\right) \right. \right. \nonumber \\
& \left. \left. \quad -\mu_j \left(2ik_x \frac{\partial \phi_j}{\partial y}-k_x^2 \Psi_j - \frac{\partial^2 \Psi_j}{\partial y^2}\right)\right]\right|_{y=y_j}=0.
\end{align}

Likewise, the continuity of displacements along the directions parallel and normal to the waveguide is given by Eqs. \eqref{continuity of transverse displacements} and \eqref{continuity of normal displacements} 

\begin{align}
\label{continuity of transverse displacements}
&\left.\left[u_x-ik_x\phi_{j}+\frac{\partial \Psi_{j}}{\partial y} \right]\right|_{y=y_j}=0, \\
\label{continuity of normal displacements}
&\left.\left[u_y-\frac{\partial \phi_{j}}{\partial y}-ik_x\Psi_{j} \right]\right|_{y=y_j}=0.
\end{align}

\section{Method} \label{Method}
\subsection{Discretisation} \label{Discretisation}

Wave solutions of the setup described in the previous section are generally split into two types - non-leaky (also known as trapped \citep{Treyssede2014FiniteWaveguides,Gallezot2017ContributionLayers,Kiefer2022ElastodynamicMetering}) and leaky modes.  In our setup, the waveguide is in contact with solids; hence the media surrounding the waveguide support shearing as well as compressive motion. For that reason, in addition to the non-leaky evanescent, wave solutions, we can further split the leaky modes solutions in two cases - those that have purely shear (shear-leaky) or shear and compressive (fully-leaky) polarisations. For these cases, the exponential growth of the amplitude of leaky waves is attributed to the exponential growth of some or all of the partial waves depicted in Fig. \ref{Schematic}. In particular, when shear-leaky waves are to be considered in an elastic half-space, the shear partial wave and corresponding transverse wave potential, $\Psi_{j}$, exponentially grow in amplitude, while the longitudinal partial wave and corresponding wave potential, $\phi_{j}$, are evanescent. In the fully-leaky case, both partial waves and their wave potentials in the solid, exhibit exponential growth. 

\begin{figure*}
    \centering
    \includegraphics[width=\textwidth]{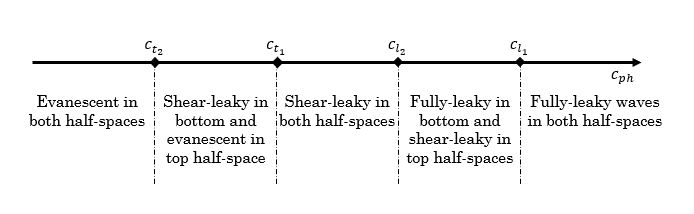}
    \caption{Illustration of the different radiation cases with an increase in phase velocity for the example case of $c_{l_1}\geq c_{l_2} > c_{t_1} \geq c_{t_2}$.}
    \label{Leaky regions}
\end{figure*}

Let us consider the structure of Fig. \ref{Schematic}, with $c_{l_1}\geq c_{l_2} > c_{t_1} \geq c_{t_2}$. In the region where the phase velocity of the propagating wave, $c_{ph}$ is less than $c_{t_2}$, evanescent waves in both solids are expected, while when it is between $c_{t_1}$ and $ c_{t_2}$, shear-leaky waves in the bottom solid and evanescent at the top are expected. When $c_{ph}$ is between $c_{l_2}$ and $ c_{t_1}$, shear-leaky waves in both spaces are expected, while when it is between $c_{l_1}$ and $c_{l_2}$, the leaky waves of the bottom solid are now expected to be fully-leaky. Finally, when $c_{ph}$ is greater than $c_{l_1}$, fully-leaky waves are expected to be radiating in both elastic half-spaces. An analogous analysis can be deduced for other combinations of longitudinal and shear bulk wavespeeds.  An illustration of this is shown in Fig. \ref{Leaky regions}. It should be noted that some exceptions of this rule, namely of subsonic radiation,  have been recorded in the literature \citep{Mozhaev2002SubsonicInterfaces,Viggen2023ModellingFluids,Kiefer2022ElastodynamicMetering,Kiefer2019CalculatingInteraction}. Such examples are still obtainable by the method to be presented, however they are not considered here. 

For the simpler case with external fluids, and a single wavespeed in the exterior, the numerical issues associated with exponential growth of the solution are discussed and a method is introduced to counter the numerical instabilities inherent to the growth \citep{Georgiades2022LeakyMethods}. As shown in Fig. \ref{Leaky regions}, with elastic half-spaces adjacent to a waveguide and all media having both shear and longitudinal wavespeeds, an exponential growth of the solution is prominent in most but not all cases, and a discretisation catering to this growth is necessary. Here we choose to solve for each of the cases described above separately, considering whether each of the partial waves of the solution exponentially grows in one or both elastic spaces. When this is the case we judiciously implement an analogous discretisation to the one described in \citep{Georgiades2022LeakyMethods} for each of Eqs. \eqref{equation of motion half-space longitudinal} and \eqref{equation of motion half-space shear}, provided that the corresponding wave potential exponentially grows. We comment that this separation into several solutions does not indicate a necessity to know the nature of the guided modes a priori; all we are doing here is composing the full solution from the several contributions of solutions for these cases. The nature of the unknown guided wave modes then determines when they emerge as solutions as we go through the calculations of the different cases.

A spectral collocation discretisation \citep{Fornberg1996AMethods,Boyd2000ChebyshevMethods}, based on Chebyshev polynomials, is applied to obtain a discrete representation of both the waveguide and the open spaces of the elastic solids. We start by utilising the Differentiation Matrix (DM) suite of Weideman and Reddy \citep{Weideman2000ASuite} to produce first and second order DMs, denoted by $D^{(1)}$ and $D^{(2)}$ respectively, for the $N$ Gauss-Lobatto Chebyshev collocation points 
\begin{equation}
    \label{chebyshev collocation points}
    s_i=-\cos\left(\frac{(i-1)\pi}{N-1}\right), 
\end{equation}

\noindent for $i=1,\dots,N$ and $s_i\in[-1,1]$. These provide an approximation of first and second order derivatives of functions defined over $[-1,1]$, evaluated at the discrete set $\{s_i\}_{i=1}^{N}$. DMs based on trigonometric or algebraic polynomials, other than Chebyshev polynomials, are also possible, with the choice of polynomials and the number of collocation points (usually points per wavelength) being crucial to the accuracy of the discretisation. Details about these choices can found in \citep{Fornberg1996AMethods,Boyd2000ChebyshevMethods,Trefethen2008SpectralMatlab}.

We require the DMs and the discretisation of the domains of definition of our displacements and wave potentials, which are $[-d,d]$ for the displacements $u_x$ and $u_y$, and the half-spaces $[-\infty,-d]$ and $[d,\infty]$ for the potentials $\phi_{1}$ and $\Psi_{1}$ and $\phi_{2}$ and $\Psi_{2}$ respectively. To achieve this, we employ mappings that relate the collocation points $\{s_i\}_{i=1}^{N}$ of $[-1,1]$ to collocation points within those domains; the corresponding DMs are then the result of the application of the chain rule.  

We start by firstly discretising the waveguide domain. A linear map, that is, a simple multiplication by $d$, maps the collocation points from $[-1,1]$ to $[-d,d]$. This provides a discretisation and DMs for the displacements in the closed domain of the waveguide. Following that, to discretise the domains of definition of the wave potentials inside the two elastic half-spaces, the rational maps

\begin{align}
    \label{h phi 1 map}
    h_{\phi_1} : s \longrightarrow y=-d-\zeta_{\phi_1}\frac{1-s}{1+s}, \\
    \label{h psi 1 map}
    h_{\Psi_1} : s \longrightarrow y=-d-\zeta_{\Psi_1}\frac{1-s}{1+s}, 
\end{align}

\noindent with $\zeta_{\phi_1},\zeta_{\Psi_1}\in \mathbb{R}^+$, map the collocation points from $[-1,1]$ to points in $[-\infty,-d]$, while the rational maps

\begin{align}
    \label{h phi 2 map}
    h_{\phi_2} : s \longrightarrow y=d+\zeta_{\phi_2}\frac{1+s}{1-s}, \\
    \label{h psi 2 map}
    h_{\Psi_2} : s \longrightarrow y=d+\zeta_{\Psi_2}\frac{1+s}{1-s}, 
\end{align}

\noindent with $\zeta_{\phi_2},\zeta_{\Psi_2}\in \mathbb{R}^+$, map the collocation points from $[-1,1]$ to points in $[d,\infty]$. 

These maps and resulting DMs, give a suitable discretisation of both elastic half-spaces, provided that the solution we are after is that of all four partial waves being evanescent in the two spaces. However, as was discussed earlier, special care and a modification of the maps of Eqs. \eqref{h phi 1 map} - \eqref{h psi 2 map}, needs to be carried out when solving for the cases of Fig. \ref{Leaky regions}, where some or all of the partial waves and corresponding wave potentials exponentially grow. When a wave potential is exponentially growing, a complex rational map is to be used instead, assigning collocation points from $[-1,1]$ to a path in the complex plane and, as will be discussed shortly, turning the exponentially growing potential into a numerically decaying function.

For clarity, we demonstrate how we modify the maps and discretise the domains with an example. Let us consider waves that radiate as shear-leaky waves in one of the half-spaces, say the bottom one, and evanescent in the other. In this case, only the transverse wave potential within the lower solid, $\Psi_{2}$, demonstrates exponential growth. As such, our attention is directed towards modifying the mapping associated with this specific potential, given in Eq. \eqref{h psi 2 map}. Allowing $\zeta_{\Psi_2}$ to take on complex values results in the discretisation of a complex path, similarly to \citep{Georgiades2022LeakyMethods}, with one of the collocation points located on the interface at $y=d$. By applying the chain rule for first and second order derivatives, we obtain the DMs for functions defined over that complex path, given in terms of the first and second order Chebyshev DMs, $D^{(1)}$ and $D^{(2)}$. Those are explicitly given by

\begin{align}
 \label{first order differentiation matrix psi 2}
    D_{\Psi_2}^{(1)}=&\text{diag}\left( \frac{2\boldsymbol{\zeta}_{\Psi_2}}{\lbrack \boldsymbol{\zeta}_{\Psi_2}+(\boldsymbol{y}_{\Psi_2}-\boldsymbol{d})\rbrack^2}\right)D^{(1)}, \\
         \label{second order differentiation matrix psi 2}
     \begin{split}
     D_{\Psi_2}^{(2)}=&\text{diag}\left(\frac{-4\boldsymbol{\zeta}_{\Psi_2}}{\lbrack \boldsymbol{\zeta}_{\Psi_2}+(\boldsymbol{y}_{\Psi_2}-\boldsymbol{d})\rbrack^3}\right)D^{(1)}\\
     &+\text{diag}\left(\frac{4\boldsymbol{\zeta}_{\Psi_2}^2}{\lbrack \boldsymbol{\zeta}_{\Psi_2}+(\boldsymbol{y}_{\Psi_2}-\boldsymbol{d})\rbrack^4}\right)D^{(2)}.
     \end{split}
\end{align}

\noindent where $\boldsymbol{\zeta}_{\Psi_2}=(\zeta_{\Psi_2},\dots,\zeta_{\Psi_2})^T$ and $\boldsymbol{d}=(d,\dots,d)^T$ are constant vectors of $N$ entries and $\boldsymbol{y}_{\Psi_2}$ is the vector of collocation points in the complex path as a result of the mapping given in Eq. \eqref{h psi 2 map}. 

The outgoing partial wave corresponding to the exponentially growing transverse wave potential $\Psi_{2}$ is characterised by the exponential term $\exp[ik_{y_{t_2}} (y-d)]$. For $\Psi_{2}$ to grow exponentially as $y$ tends to $\infty$, $-\Im(k_{y_{t_2}})$ must be a positive real number. By performing the discretisation presented here, we are mapping $y$ to a complex path such that the exponent may be expressed as $ik_{y_{t_2}} (y'-d)$, with $y'=\zeta_{\Psi_{2}}(y-d)+d$. With an appropriate choice of $\zeta_{\Psi_{2}}$ \citep{Georgiades2022LeakyMethods}, the exponential term can then be made to decay. 

This newly constructed, numerically decaying, function is the analytical continuation of the transverse potential $\Psi_{2}$ over the complex path defined by Eq. \eqref{h psi 2 map}, and we denote it by $\tilde {\Psi}_{2}$. Just like $\Psi_{2}$, it is fully determined by its value at the interface, while due to our choice of mapping, which maintains one collocation point at the physical interface of the waveguide and the bottom solid, $\Psi_{2}$ and $\tilde {\Psi}_{2}$ attain the same value at that interface. What is more, as a time harmonic function, $\tilde {\Psi}_{2}$ also satisfies the governing Eq. \eqref{equation of motion half-space shear} for the same wavenumber $k_x$. As a result, by discretising a complex path as opposed to the physical half-space of the elastic solid, we are able to solve Eq. \eqref{equation of motion half-space shear} for the numerically stable and decaying function, $\tilde {\Psi}_{2}$, and the same complex wavenumber $k_x$, and then reconstruct the physical, exponentially growing, solution of $\Psi_{2}$ from its value at the interface. We note the resemblance of the complex rational maps of Eqs. \eqref{h phi 1 map} - \eqref{h psi 2 map} to the complex coordinate stretching used in the technique of PMLs; see the Appendix for a discussion.

By employing this technique, Eq. \eqref{equation of motion half-space shear} can be discretised and presented in matrix form as follows 
\begin{equation}
    \label{equation of motion half-space shear discretised}
    \left( \rho_{j}\omega^2 I-k_x^2\mu_{j}I+\mu_{j}D_{\Psi_2}^{(2)} \right) \boldsymbol{\tilde{\Psi}_{2}} =\boldsymbol{0}, 
\end{equation}
where $I$ is the $N\times N$ identity matrix and $\boldsymbol{\tilde{\Psi}_{2}}$ is the vector containing the evaluations of $\tilde{\Psi}_{2}$ at the collocation points in the complex path. 

For the rest of the evanescent partial waves and their decaying wave potentials, we use the rational maps of Eqs. \eqref{h phi 1 map}-\eqref{h psi 1 map}. This time, real valued parameters are used for the maps, as there is no need to impose further decay to the already decaying potentials we seek, and analogous discretisations of the domains and matrix formulations of Eqs. \eqref{equation of motion waveguide x}-\eqref{equation of motion half-space shear} are obtained.

Finally, similarly to \cite{Georgiades2022LeakyMethods, Adamou2004SpectralMedia, Kiefer2019CalculatingInteraction}, we merge all discretised equations into a single matrix equation, replacing the rows corresponding to the collocation points at the two interfaces with the relevant, discretised, interface conditions, and those at infinity with Dirichlet zero boundary conditions. This results in a polynomial eigenvalue problem
\begin{equation}
    (k_x^2L_2+k_xL_1+L_0)\boldsymbol{u}=0,
    \label{eigenvalue problem}
\end{equation}
where $\boldsymbol{u}=(\boldsymbol{\phi_{1}},\boldsymbol{\Psi_{1}},\boldsymbol{u_x},\boldsymbol{u_y},\boldsymbol{\phi_{2}},\boldsymbol{\tilde{\Psi}_{2}})^T$, with each entry in the vector $\boldsymbol{u}$ being the vector of the respective quantity evaluated at the collocation points, and where $L_2, L_1$ and $L_0$ are the $k_x^2$-dependent, $k_x$-dependent and $k_x$-independent components of the merged matrix equation. Due to the choice of complex transformation, the polynomial eigenvalue problem of Eq. \eqref{eigenvalue problem} only admits pairs of positive and negative phase velocity outgoing wave solutions while, for completeness, the corresponding incoming wave solutions are also attainable via the same method with a choice of the conjugate complex parameters.

We treat other cases of Fig. \ref{Leaky regions} in the same manner, each time discretising a path in the complex plane to accommodate for the exponential growth of leaky waves. We solve each of the resulting eigenvalue problems using companion matrix linearisations \citep{Kiefer2019CalculatingInteraction,Mackey2006StructuredLinearizations} and conventional eigenvalue problem solvers. 

\subsection{Finite Element Modelling}

As an additional verification  tool, we have used FE modelling to model wave propagation in the waveguide of finite thickness adjacent to two elastic half-spaces shown schematically in Fig. \ref{Schematic}. 

We make use of the GPU-based FE solver Pogo\cite{Huthwaite2014AcceleratedGPU}, that discretises the equation of motion from elasticity theory in space and time. Discretisation in time is achieved by the finite-difference method\citep{BatheK.J.1996FiniteProcedures} and in this way displacement values at a given time step are retrieved as a function of displacement values at the two time steps preceding it. The FE modelling approach described here follows the FE models applied in other elastic wave progation studies \citep{Sarris2021AttenuationRoughness,Sarris2023UltrasonicDamage,Haslinger2019AppraisingDefects,Huang2020MaximizingPolycrystals}.

\begin{figure}
    \centering
    \includegraphics[width=\columnwidth]{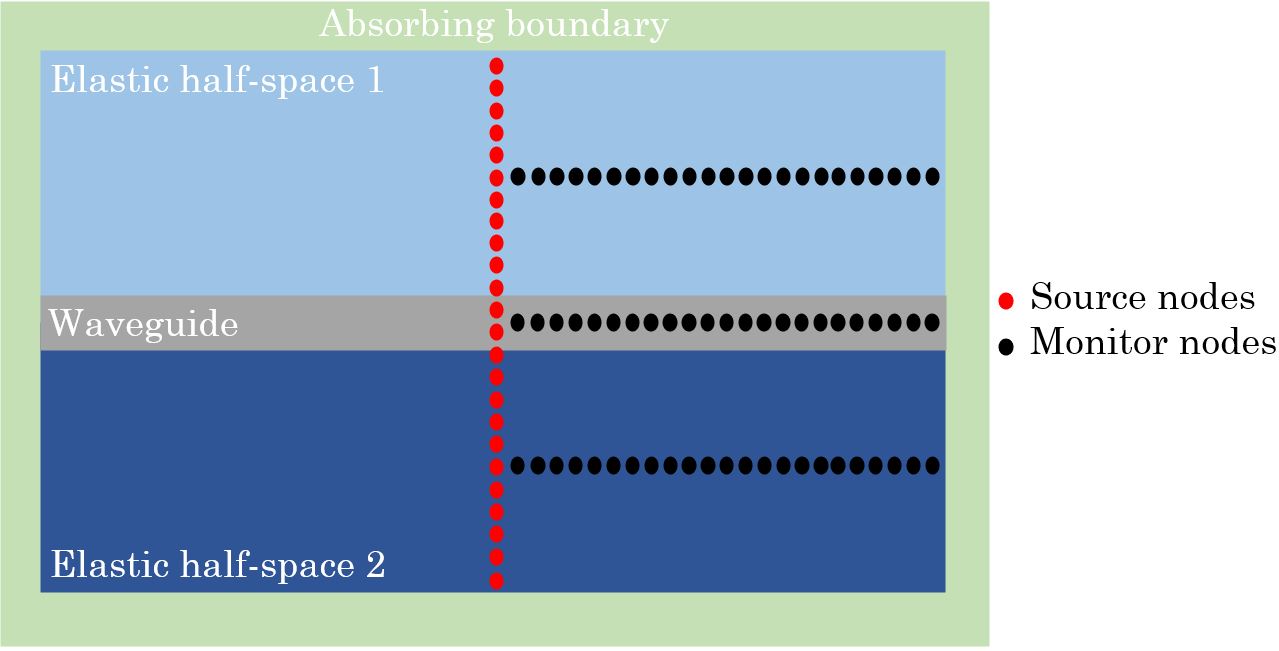}
    \caption{(Colour online) A schematic of a FE domain used to model elastic wave propagation and retrieve the dispersion curves of guided modes in a waveguide adjacent to two elastic half-spaces. In the figure, source nodes are shown in red and monitor nodes in black. Absorbing layers around the domain are shown in green while different colours denote the different materials within the domain.}
    \label{Scehmatic FE}
\end{figure}

For the FE modelling of leaky wave propagation in a waveguide coupled to two elastic half-spaces, a rectangular domain was used - a schematic of the cross sectional area of the domain is shown in Fig. \ref{Scehmatic FE}. Each element of the domain is assigned with the material properties of the space it is a part of (Young's modulus, density, Poisson's ratio), while linear, square elements of four nodes each were used. To ensure sufficient convergence of the solution \citep{Drozdz2008EfficientMedia}, the mesh element size was set to approximately $dx=\lambda/30$. Here, $\lambda$ is the smallest bulk wavelength amongst the three materials, calculated at the frequency of the simulation. Stiffness-reducing absorbing layers \citep{Rajagopal2012OnPackages} were defined around the domain to accommodate for its unbounded nature, with thickness equal to four times the longest bulk wavelength. In the absorbing layers region, damping proportional to the stiffness is gradually introduced into the system throughout the thickness of the region, resulting in the attenuation of outgoing waves and minimal reflections.

A source line, comprised of multiple source nodes, was defined in the middle of the domain - source nodes are nodes with predefined displacement/force, and are used to excite the propagating wave in the model. Two Hann-windowed tonebursts were applied to each source node, with each toneburst displacing the node in the $x$ or the $y$ direction. Depending on the nature of the desired excitation, multi-modal or single-mode excitation, the amplitude and relative phase of the tonebursts at each source node were determined. For a multi-modal excitation, all source nodes were excited in phase with the same broadband (single-cycled) toneburst in both directions.
For single-mode excitation, the method discussed in \cite{Pavlakovic1998LeakyNDT} was followed; first, the displacement fields (mode shapes) of the desired mode were extracted from our method and verified with those from DISPERSE. Then, the amplitude and phase of the mode shapes were used to modulate the amplitudes and phases of the excitation signals accordingly. It is worth noting that due to the use of a truncated domain and source node placement within the domain, but excluding absorbing layers, the proposed excitation violates the orthogonality relations of eigenfunctions of the system. This results in a non-pure excitation, possibly exciting other modes as well as the desired one. Since our primary objective was the qualitative validation of SCM results, achieving a pure mode excitation was not pursued. For an analysis of the modal basis of such systems and the orthogonality of eigenfunctions for pure-mode excitations, see \citep{Gallezot2018AWaveguides,Sammut1976LeakyExcitation,Snyder1983OpticalTheory,Kiefer2022ElastodynamicMetering} and references therein.

A monitor line, comprised of multiple monitor nodes, was placed horizontally, to the right of the source line, along the length of the domain, as shown in Fig. \ref{Scehmatic FE} - a monitor node is a node whose displacement profile is recorded and stored for post-processing. The placement of the monitor line is dependent on the nature of the desired mode to be observed. Leaky modes attenuate as they propagate along the waveguide, with most of their energy observable in the exterior solids, requiring a monitor line that is placed close to the source line and in the exterior elastic solids. Contrary to that, non-leaky modes remain trapped in the waveguide and can propagate longer distances; thus a monitor line in the middle of the waveguide is more suitable for their observation. In post-processing, a two-dimensional fast Fourier transform (2D-FFT)\citep{Alleyne1991ASignals}, in space and time, is applied to the collection of the stored monitor signals of interest. The results of the 2D-FFT were then normalised to the maximum amplitude of all the modes recorded to obtain the dispersion curves of the excited modes in a normalised amplitude scale.

A typical model size is of the order of $1\times 10^7$ degrees of freedom. The models were solved using an Nvidia (Santa Clara, California) 2080 Ti graphics card with 8GB of memory. The solution time was approximately forty seconds. More details about the excitation, the size of the domain and the material properties will be given in the following section as they are specific to the example presented there. 

\section{Results and Discussion} \label{Results and Discussion}
We now utilise our SCM method to retrieve the dispersion and attenuation curves of leaky and non-leaky modes in an elastic waveguide in contact with two elastic half-spaces. The obtained dispersion and attenuation plots are then verified against the commercially-available software DISPERSE while the dispersion curves are also validated against FE simulations of elastic wave propagation using Pogo. 

Without loss of generality, we illustrate the method on the example case of adhesive joints of two aluminium alloys, studied in \citep{Lowe1994TheJoints} to demonstrate ways of determining the material properties and thickness of the adhesive layer and assess the quality of the joint; an example inspired by real NDE applications. An epoxy adhesive is modelled as a waveguide of finite half-thickness $d=0.5$ mm, density $\rho_{\text{ep}}=1.17 \text{ g/cm}^3$, longitudinal wavespeed $c_{l_{\text{ep}}}=2.61 \text{ m/ms}$ and transverse wavespeed $c_{t_{\text{ep}}}=1.1 \text{ m/ms}$. The spaces of the aluminium components are modelled as linear elastic half-spaces of isotropic media, both with density $\rho_{\text{al}}=2.82 \text{ g/cm}^3$, longitudinal wavespeed $c_{l_{\text{al}}}=6.33 \text{ m/ms}$ and transverse wavespeed $c_{t_{\text{al}}}=3.12 \text{ m/ms}$.

For the retrieval of the solution, the phase velocity domain was split up in three regions in accordance with Fig. \ref{Leaky regions} and the wavespeeds of the adjacent aluminium half-spaces. For each region, an angular frequency sweep of 150 frequency steps, up to $\omega=30 \text{ MHz}$, was performed, with the complex parameter pairs $(\zeta_{\phi_j},\zeta_{\Psi_j})$ being chosen to accommodate the radiation characteristics of the solutions to be found within that particular region and remaining constant for all frequencies of the sweep. In line with the example of shear-leaky waves discussed in the previous section, the pairs are appropriately chosen for each region to distribute the collocation points across real and complex paths and ensure the numerical decay of all the wave potentials to be found. For the retrieval of non-leaky solutions, $\left(\zeta_{\phi_j},\zeta_{\Psi_j}\right)=\left(10,10\right)$ was chosen, while $\left(\zeta_{\phi_j},\zeta_{\Psi_j}\right)=\left(10,10i\right)$ were chosen for shear-leaky and $\left(\zeta_{\phi_j},\zeta_{\Psi_j}\right)=\left(10i,10i\right)$ for fully-leaky ones. In each domain, $50$ collocation points were used. The results of our SCM were found to be stable to the choices of complex parameters and number of collocation points.  

Subsequently, the numerical solutions obtained underwent a comprehensive filtering process; only modes of finite attenuation and positive phase velocity were kept. In terms of the wavenumber, $k_x$, that is achieved by requiring $\Re(k_x)>0$, which ensures a mode of non-negative and finite phase velocity, and $|k_x|<\infty$, which results in non-zero phase velocity and finite attenuation modes. We further filter out highly attenuating modes, by requiring a non-negative attenuation of less than $15$ Np/mm. Lastly, we get rid of spurious, non-converged and bulk wave solutions by requiring that the mode complies with interface conditions of Eqs. \eqref{continuity of normal stress}-\eqref{continuity of normal displacements} (up to 3 decimal places), under the assumption that the wave potentials of the outgoing waves in the two elastic half-spaces assume the general forms $\phi_1=A_1\exp\left[-ik_{y_{l_1}}(y+d)\right]$, $\Psi_1=B_1\exp\left[-ik_{y_{t_1}}(y+d)\right]$, $\phi_2=A_2\exp\left[ik_{y_{l_2}}(y-d)\right]$ and $\Psi_2=B_2\exp\left[ik_{y_{t_2}}(y-d)\right]$, subject to the choice of square root of $k_{y_{l_j}}$ and $k_{y_{t_j}}$, for which $\tilde {\phi}_{j}$ and $\tilde {\Psi}_{j}$ correspond to the solution found by the SCM. The computation took a total of approximately 18 minutes on a personal computer with 32GB RAM memory and an i7-10700 CPU.

\begin{figure}
    \centering
    \includegraphics[width=\columnwidth]{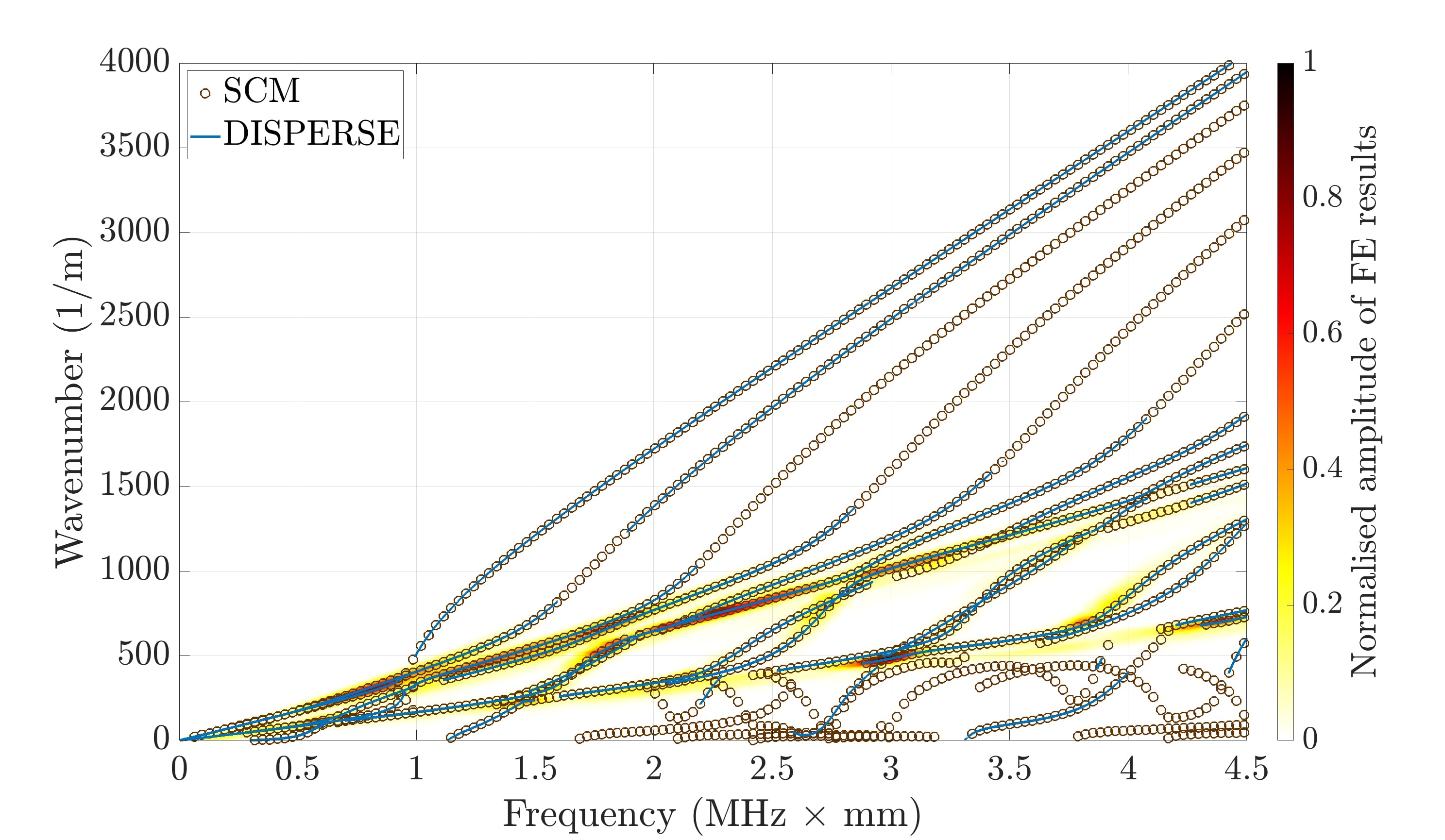}
    \caption{(Colour online) Dispersion curves of an epoxy waveguide in contact with an aluminium half-space on each side, obtained by our SCM (denoted by $\circ$), by DISPERSE (denoted by --) and by FE modelling (normalised intensity in colour). }
    \label{Multi-modal excitation results}
\end{figure}

A comparison of the SCM results with those from a FE simulation of a multi-modal excitation of our adhesive joint domain was possible from the excitation of all source nodes in phase and with the same broadband toneburst in both the $x$ and $y$ directions. The excitation toneburst used at each source node was the sum of sixteen single-cycled Hann-windowed toneburtsts, each with a central frequency varying from $0.25-4$ MHz at $0.25$ MHz increments. In this case, wavelengths and element size, $dx$, were calculated at $4$ MHz. Employing a toneburst excitation composed of the cumulative energy from 16 distinct broadband tonebursts is equivalent to exciting each of these 16 tonebursts independently and subsequently combining their excited fields. This ensures that at each frequency increment, modes are excited with more intensity than they would have been excited with, say, a single-cycled toneburst centred at the middle of the frequency range of interest, thereby promoting stronger mode excitation throughout the desired frequency range. The domain, excluding absorbing layers, was 11 mm thick, 1 mm of epoxy in the middle of the domain and 5 mm of aluminium each side of the waveguide, and 100 mm in length. The monitor line used for post-processing was placed at $y=2$ mm. Normalised results of the described FE model are shown in Fig. \ref{Multi-modal excitation results}.

As shown in Fig. \ref{Multi-modal excitation results}, there is good agreement between all three methods for the modes that the chosen forcing in the FE model was able to excite. There were, however, instances where DISPERSE was not able to fully retrieve a mode, or missed it altogether, highlighting the merits of using the proposed SCM instead of root finding techniques. 

\begin{figure} [H]
    \centering
    \includegraphics[width=\columnwidth]{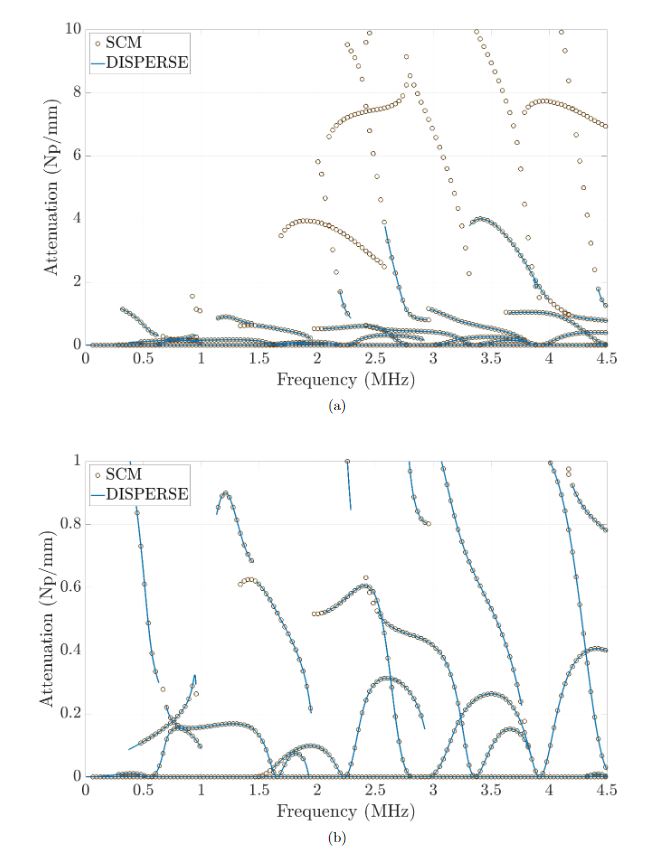}
    \caption{(Colour online) Comparison of attenuation values obtained by our SCM (denoted by $\circ$) and by DISPERSE (denoted by --) for waves in an epoxy waveguide in contact with an aluminium half-space on each side \hyperref[Attenuation plot]{\ref{Attenuation plot}a}. A detail of the region with attenuation of up to 1 Np/mm \hyperref[Attenuation plot]{\ref{Attenuation plot}b}.}
    \label{Attenuation plot}
\end{figure}

Additionally, a comparison of the attenuation values obtained by the SCM against those of DISPERSE is shown in Fig. \hyperref[Attenuation plot]{\ref{Attenuation plot}a}. For visualisation purposes, solutions with high attenuation, $\Im(k_x) > 10 \text{ Np/mm}$, are not displayed here. A detail of the region of up to 1 Np/mm is shown in Fig. \hyperref[Attenuation plot]{\ref{Attenuation plot}b}. The attenuation values obtained by our SCM match those obtained by DISPERSE, demonstrating in this way the ability of our SCM to accurately capture the complex wavenumber solutions of the leaky modes of the system.

Following these initial validations of the SCM results, we proceed with presenting two cases of single-mode excitation - one where both DISPERSE and SCM were able to retrieve the mode to be excited, and one where only SCM was able to obtain the relevant mode. In each case, the excitation tonebursts were set to the same central frequency, and for the purposes of a narrowband excitation, the number of cycles in the excitation signal was increased to eleven. The size of the absorbing layers, and $dx$ were calculated accordingly at the central frequency. 

\begin{figure} [H]
    \centering
    \includegraphics[width=\columnwidth]{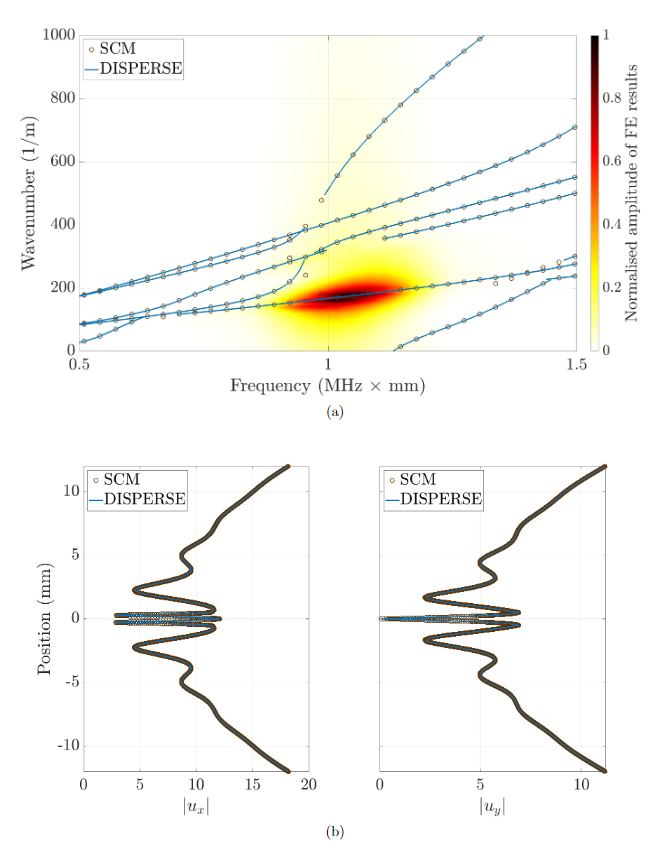}
\caption{\label{single-mode excitation scm and disperse} (Colour online) Dispersion curves showing the normalised intensity and dispersion of the shear-leaky mode excited at $f=1.03$ MHz and $c_{ph}=5.98$ m/ms with a single-mode excitation, validated against SCM and DISPERSE dispersion curves \hyperref[single-mode excitation scm and disperse]{\ref{single-mode excitation scm and disperse}a}. The exponentially growing mode shapes from the SCM that were used for the excitation are validated against those obtained by DISPERSE \hyperref[single-mode excitation scm and disperse]{\ref{single-mode excitation scm and disperse}b}.  }
\end{figure}

Fig. \hyperref[single-mode excitation scm and disperse]{\ref{single-mode excitation scm and disperse}a} shows the results from the independent excitation of a leaky mode using FE modelling, at $f=1.03$ MHz and $c_{ph}=5.98$ m/ms, which both DISPERSE and SCM were able to retrieve. The FE results shown are from the post-processing of monitor signals from a monitor line placed horizontally at $y=5$ mm. As the phase velocity of the mode is between the transverse and longitudinal bulk wavespeeds of the exterior aluminium spaces, from Fig. \ref{Leaky regions}, this is in fact a shear-leaky mode, radiating energy in both exterior domains purely in the form of shear waves. To achieve this excitation and validate our method, mode shapes from the SCM were used and are displayed in Fig. \hyperref[single-mode excitation scm and disperse]{\ref{single-mode excitation scm and disperse}b}, where they can be seen to exponentially grow away from the waveguide. Alongside them, mode shapes from DISPERSE are also shown, demonstrating a pleasing concordance between the two methods in computing the mode shapes of the shear-leaky mode.

\begin{figure} [H]
    \centering
    \includegraphics[width=\columnwidth]{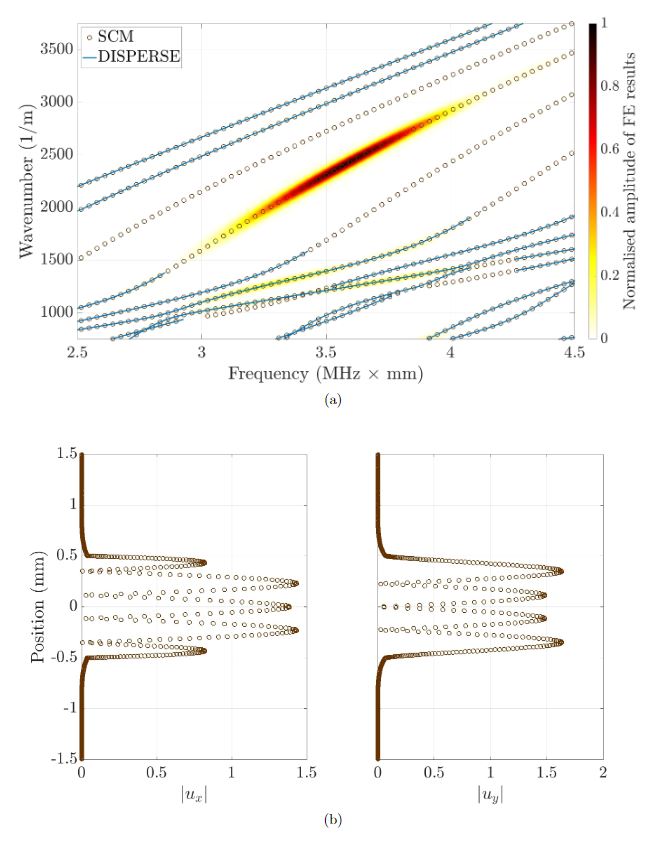}
\caption{\label{single-mode excitation scm only} (Colour online) Dispersion curves showing the normalised intensity and dispersion of the non-leaky mode, found only by SCM, and excited at $f=3.53$ MHz and $c_{ph}=1.50$ m/ms with a single-mode excitation \hyperref[single-mode excitation scm only]{\ref{single-mode excitation scm only}a}. The decaying mode shapes from the SCM that were used for the excitation are shown in Fig. \hyperref[single-mode excitation scm only]{\ref{single-mode excitation scm only}b}.}
\end{figure}

Fig. \hyperref[single-mode excitation scm only]{\ref{single-mode excitation scm only}a} shows the results from the excitation of a non-leaky mode, at $f=3.53$ MHz and $c_{ph}=1.50$ m/ms, which only our SCM was able to retrieve, alongside the relevant FE result. As the energy of such modes is confined within the boundaries of the waveguide, the monitor line of the FE model is now placed horizontally in the middle of the waveguide. The mode shapes used for this excitation, decaying away from the waveguide, are shown in Fig. \hyperref[single-mode excitation scm only]{\ref{single-mode excitation scm only}b} as obtained from our SCM.

It is important to note that any input to the FE model will inherently result in a weighted superposition of its modal responses \cite{BatheK.J.1996FiniteProcedures}. Consequently, we regard the agreement between the FE model and the results obtained from SCM and DISPERSE in the figures above as a robust validation of both methods. This validation, in particular, emphasises the reliability of SCM in scenarios where a DISPERSE solution is not readily available. 

\newpage

\section{Conclusion}\label{Conclusion}
Leaky Lamb wave propagation in a waveguide that is in contact with elastic half-spaces has been considered; these wave modes are often characterised by a complex wavenumber, a subsequent exponential growth in amplitude far into the surrounding elastic media and a highly attenuative nature along the waveguide, making their retrieval by traditional numerical solvers challenging. Here, we overcome these obstacles by formulating a theoretical framework for the problem, discretising it using a SCM, that also accommodates for the exponentially growing nature of the problem, and finally solving the resulting polynomial eigenvalue problem for the retrieval of the wave mode solutions of interest. The framework introduced was developed for the general case of a waveguide between two possibly different elastic half-spaces, however it can be simplified to handle other cases of loading by elastic half-spaces or inviscid fluids, like waveguides loaded only on one side by a solid or one with a solid and a fluid half-space each side of the waveguide.

The results of our SCM have been verified against a commercially available software, DISPERSE, and FE modelling in places where a DISPERSE solution was not readily available. The SCM approach is accurate, versatile, it readily identifies all leaky modes in elasticity without requiring knowledge of the modes or any tuning a priori and is easy to code; the key point is to remove the exponential growth from the numerical scheme, whilst it remains within the physics, and doing so via the mapping to a complex path reduces the leaky problem to a standard form for numerical implementation.

\begin{acknowledgments}
The authors would like to thank the reviewers for their valuable contributions in improving the quality of this work. The research of EG and RVC has received funding from the European Union's Horizon 2020 FET Open programme under grant agreement No. 863179 (BOHEME).  
\end{acknowledgments}

\section{Author Declaration}
The authors have no conflicts of interest
to disclose. 
\section{Data Availability}
The data that support the findings of this study are available from the corresponding author upon reasonable request.

\appendix*
\section{Relationship to PMLs} 

PMLs were originally defined by Berenger \citep{Berenger1994AWaves} for the simulation of free-space for electromagnetic scattering problems using the finite-difference time-domain method. Therein, a layer of perfectly matched impedance was placed around the edges of the finite computational domain. The layer is defined through a change in the magnetic and electric conductivities; these were chosen to have a power law dependence on the ratio of distance from the interface between PML and unaltered material to the thickness of the PML: 

\begin{equation} \label{Berenger}
    \sigma(r)=\sigma_m \left(\frac{r}{\delta}\right)^n,
\end{equation}

\noindent with $\sigma_m$ the limiting magnetic or electric conductivity inside the PML, $r$ the distance from the interface inside the matched layer, $\delta$ the PML's thickness and $n$ a non-negative integer. The choice of conductivities in Eq. \eqref{Berenger} ensure that waves enter the PML without reflections and are absorbed; in practical terms the choices of the parameters $\sigma_m,\delta$ and $n$ are all important in achieving the optimal outcome.

Later, Chew and Weedon \citep{Chew1994ACoordinates} proposed an alternative viewpoint of PML using the modified Maxwell's equations in terms of a complex coordinate stretching. The same idea was extended to elastodynamics \citep{Chew1996PerfectlyCondition,Basu2003PerfectlyImplementation,Kim2009TheLayer} to produce a medium that absorbs outgoing waves. The coordinate stretching is generally (see for instance \citep{Treyssede2014FiniteWaveguides,Kim2009TheLayer,Chew1994ACoordinates,Collino2001ApplicationMedia,Harari2006StudiesWaves}) given by 

\begin{equation} \label{stretching}
    r\rightarrow \tilde{r}= \int_{0}^{r} \gamma (s) \text{ d}s,
\end{equation}
with $\gamma$ a complex-valued function such that

\begin{equation} \label{conditions on gamma}
    \begin{cases}
        \Im\{\gamma(r)\}>0 & \text{inside the PML,} \\
        \gamma(r)=1 & \text{elsewhere.}
    \end{cases}
\end{equation}

\noindent A typical choice for $\gamma$ in Eq. \eqref{stretching} that satisfies the conditions of Eq. \eqref{conditions on gamma} is 

\begin{equation} \label{gamma}
    \gamma (r) = 1+\frac{i}{\omega}\sigma(r),
\end{equation}

\noindent with $\sigma(r)$ taken from Berenger's conductivity Eq. \eqref{Berenger} inside the PML and zero everywhere else. In essence the undisturbed material is surrounded by an impedance matched layer with material properties gradually changing throughout the layer to minimise any reflections back into the computational domain of interest; constant, rather than gradually changing, choices, i.e. $n=0$, tend to perform poorly. 

In theory PMLs are layers of perfectly matched impedance, in practice, reflections can occur from the interface between the PML and unaltered material due to the discretisation of the wave equation for use in a computer simulation \citep{Rajagopal2012OnPackages,Pled2022ReviewDomains}. Skelton \textit{et al.}\citep{Skelton2007GuidedLayers} provide the mathematical grounds for those reflections, where they can be interpreted as emitted disturbances propagating away from a collection of sources placed at the elastic-PML interface; it was concluded that to minimise reflections, the power dependence of Eq. \eqref{Berenger} must at least be a quadratic ($n \geq 2)$. 

Turning away from diffraction and scattering to now consider eigenmodes of structure, 
Treyssède \textit{et al.}\citep{Treyssede2014FiniteWaveguides} used PMLs surrounding a region of unaltered  material for the computation of modes of an open stratified waveguide; results from a constant PML $(n=0)$ - previously adopted in the literature for open waveguides \citep{Pelat2011AWaveguides,Bonnet-BenDhia2010AWaveguides,Agha2008OnFibres} - and from a parabolic PML $(n=3)$ were compared and computational accuracy improved when the continuously increasing absorbing function of the parabolic PML was used. 

The methodology used in the body of this article can be viewed as a constant PML (or of constant conductivity under Berenger's definition), like that in \citep{Treyssede2014FiniteWaveguides}, but now attached directly onto the waveguide. Considering $n=0$ in Eq. \eqref{Berenger}, the complex coordinate stretching of Eq. \eqref{stretching}, for a PML attached at the waveguide interface at $y=d$, becomes 
\begin{equation} \label{constant stretching}
    r \rightarrow \tilde{r}=d+\left(1+\frac{i}{\omega}\sigma_m \right)r,
\end{equation}
inside the PML. With a choice of 
\begin{equation} \label{sigma_m}
    \sigma_m=-i\omega \left(\zeta_{\Psi_2}-1\right), 
\end{equation}
and $r=y-d$, Eq. \eqref{constant stretching} maps collocation points from the real line $[d,\infty]$ to the complex path defined by Eq. \eqref{h psi 2 map}. Crucially, now both the elastic medium and the PML have collocation points on the same physical space of the waveguide-exterior interface and, by enforcing continuity of tractions and displacements, there is no impedance mismatch at the elastic-PML interface; the disadvantages of a constant PML are negated as reflections from the elastic-PML interface do not pollute the computation and the constant map effectively computes leaky waves. 



\bibliography{references}

\begin{thebibliography}{68}%
\makeatletter
\providecommand \@ifxundefined [1]{%
 \@ifx{#1\undefined}
}%
\providecommand \@ifnum [1]{%
 \ifnum #1\expandafter \@firstoftwo
 \else \expandafter \@secondoftwo
 \fi
}%
\providecommand \@ifx [1]{%
 \ifx #1\expandafter \@firstoftwo
 \else \expandafter \@secondoftwo
 \fi
}%
\providecommand \natexlab [1]{#1}%
\providecommand \enquote  [1]{``#1''}%
\providecommand \bibnamefont  [1]{#1}%
\providecommand \bibfnamefont [1]{#1}%
\providecommand \citenamefont [1]{#1}%
\providecommand \href@noop [0]{\@secondoftwo}%
\providecommand \href [0]{\begingroup \@sanitize@url \@href}%
\providecommand \@href[1]{\@@startlink{#1}\@@href}%
\providecommand \@@href[1]{\endgroup#1\@@endlink}%
\providecommand \@sanitize@url [0]{\catcode `\\12\catcode `\$12\catcode `\&12\catcode `\#12\catcode `\^12\catcode `\_12\catcode `\%12\relax}%
\providecommand \@@startlink[1]{}%
\providecommand \@@endlink[0]{}%
\providecommand \url  [0]{\begingroup\@sanitize@url \@url }%
\providecommand \@url [1]{\endgroup\@href {#1}{\urlprefix }}%
\providecommand \urlprefix  [0]{URL }%
\providecommand \Eprint [0]{\href }%
\providecommand \doibase [0]{https://doi.org/}%
\providecommand \selectlanguage [0]{\@gobble}%
\providecommand \bibinfo  [0]{\@secondoftwo}%
\providecommand \bibfield  [0]{\@secondoftwo}%
\providecommand \translation [1]{[#1]}%
\providecommand \BibitemOpen [0]{}%
\providecommand \bibitemStop [0]{}%
\providecommand \bibitemNoStop [0]{.\EOS\space}%
\providecommand \EOS [0]{\spacefactor3000\relax}%
\providecommand \BibitemShut  [1]{\csname bibitem#1\endcsname}%
\let\auto@bib@innerbib\@empty
\bibitem [{\citenamefont {Shen}\ \emph {et~al.}(2014)\citenamefont {Shen}, \citenamefont {Hirose},\ and\ \citenamefont {Yamaguchi}}]{Shen2014DispersionStudy}%
  \BibitemOpen
  \bibfield  {author} {\bibinfo {author} {\bibfnamefont {Y.}~\bibnamefont {Shen}}, \bibinfo {author} {\bibfnamefont {S.}~\bibnamefont {Hirose}},\ and\ \bibinfo {author} {\bibfnamefont {Y.}~\bibnamefont {Yamaguchi}},\ }\bibfield  {title} {\bibinfo {title} {{Dispersion of ultrasonic surface waves in a steel-epoxy-concrete bonding layered medium based on analytical, experimental, and numerical study}},\ }\href {https://doi.org/10.1016/j.csndt.2014.07.002} {\bibfield  {journal} {\bibinfo  {journal} {Case Stud. Nondestruct. Test. Evaluation}\ }\textbf {\bibinfo {volume} {2}},\ \bibinfo {pages} {49} (\bibinfo {year} {2014})}\BibitemShut {NoStop}%
\bibitem [{\citenamefont {Lellouch}\ \emph {et~al.}(2021)\citenamefont {Lellouch}, \citenamefont {Biondi}, \citenamefont {Biondi}, \citenamefont {Luo}, \citenamefont {Jin},\ and\ \citenamefont {Meadows}}]{Lellouch2021PropertiesFiber}%
  \BibitemOpen
  \bibfield  {author} {\bibinfo {author} {\bibfnamefont {A.}~\bibnamefont {Lellouch}}, \bibinfo {author} {\bibfnamefont {E.}~\bibnamefont {Biondi}}, \bibinfo {author} {\bibfnamefont {B.~L.}\ \bibnamefont {Biondi}}, \bibinfo {author} {\bibfnamefont {B.}~\bibnamefont {Luo}}, \bibinfo {author} {\bibfnamefont {G.}~\bibnamefont {Jin}},\ and\ \bibinfo {author} {\bibfnamefont {M.~A.}\ \bibnamefont {Meadows}},\ }\bibfield  {title} {\bibinfo {title} {{Properties of a deep seismic waveguide measured with an optical fiber}},\ }\bibfield  {journal} {\bibinfo  {journal} {Phys. Rev. Res.}\ }\textbf {\bibinfo {volume} {3}},\ \href {https://doi.org/10.1103/PhysRevResearch.3.013164} {10.1103/PhysRevResearch.3.013164} (\bibinfo {year} {2021})\BibitemShut {NoStop}%
\bibitem [{\citenamefont {Chimenti}\ and\ \citenamefont {Nayfeh}(1985)}]{Chimenti1985LeakyLaminates}%
  \BibitemOpen
  \bibfield  {author} {\bibinfo {author} {\bibfnamefont {D.~E.}\ \bibnamefont {Chimenti}}\ and\ \bibinfo {author} {\bibfnamefont {A.~H.}\ \bibnamefont {Nayfeh}},\ }\bibfield  {title} {\bibinfo {title} {{Leaky Lamb waves in fibrous composite laminates}},\ }\href {https://doi.org/10.1063/1.336268} {\bibfield  {journal} {\bibinfo  {journal} {J. Appl. Phys.}\ }\textbf {\bibinfo {volume} {58}},\ \bibinfo {pages} {4531} (\bibinfo {year} {1985})}\BibitemShut {NoStop}%
\bibitem [{\citenamefont {Nagy}(1994)}]{Nagy1994LeakyMaterials}%
  \BibitemOpen
  \bibfield  {author} {\bibinfo {author} {\bibfnamefont {P.~B.}\ \bibnamefont {Nagy}},\ }\bibfield  {title} {\bibinfo {title} {{Leaky Guided Wave Propagation Along Imperfectly Bonded Fibers in Composite Materials}},\ }\bibfield  {journal} {\bibinfo  {journal} {J. Nondestruct. Eval.}\ }\textbf {\bibinfo {volume} {13}},\ \href {https://doi.org/10.1007/BF00728251} {10.1007/BF00728251} (\bibinfo {year} {1994})\BibitemShut {NoStop}%
\bibitem [{\citenamefont {Long}\ \emph {et~al.}(2003)\citenamefont {Long}, \citenamefont {Lowe},\ and\ \citenamefont {Cawley}}]{Long2003AttenuationPipes}%
  \BibitemOpen
  \bibfield  {author} {\bibinfo {author} {\bibfnamefont {R.}~\bibnamefont {Long}}, \bibinfo {author} {\bibfnamefont {M.~J.~S.}\ \bibnamefont {Lowe}},\ and\ \bibinfo {author} {\bibfnamefont {P.}~\bibnamefont {Cawley}},\ }\bibfield  {title} {\bibinfo {title} {{Attenuation characteristics of the fundamental modes that propagate in buried iron water pipes}},\ }\href {https://doi.org/10.1016/S0041-624X(03)00166-5} {\bibfield  {journal} {\bibinfo  {journal} {Ultrasonics}\ }\textbf {\bibinfo {volume} {41}},\ \bibinfo {pages} {509} (\bibinfo {year} {2003})}\BibitemShut {NoStop}%
\bibitem [{\citenamefont {Long}\ \emph {et~al.}(2004)\citenamefont {Long}, \citenamefont {Vogt}, \citenamefont {Lowe},\ and\ \citenamefont {Cawley}}]{Long2004MeasurementWaveguide}%
  \BibitemOpen
  \bibfield  {author} {\bibinfo {author} {\bibfnamefont {R.}~\bibnamefont {Long}}, \bibinfo {author} {\bibfnamefont {T.}~\bibnamefont {Vogt}}, \bibinfo {author} {\bibfnamefont {M.~J.~S.}\ \bibnamefont {Lowe}},\ and\ \bibinfo {author} {\bibfnamefont {P.}~\bibnamefont {Cawley}},\ }\bibfield  {title} {\bibinfo {title} {{Measurement of acoustic properties of near-surface soils using an ultrasonic waveguide}},\ }\href {https://doi.org/10.1190/1.1707065} {\bibfield  {journal} {\bibinfo  {journal} {Geophys.}\ }\textbf {\bibinfo {volume} {69}},\ \bibinfo {pages} {460} (\bibinfo {year} {2004})}\BibitemShut {NoStop}%
\bibitem [{\citenamefont {Leinov}\ \emph {et~al.}(2016)\citenamefont {Leinov}, \citenamefont {Lowe},\ and\ \citenamefont {Cawley}}]{Leinov2016UltrasonicPipes}%
  \BibitemOpen
  \bibfield  {author} {\bibinfo {author} {\bibfnamefont {E.}~\bibnamefont {Leinov}}, \bibinfo {author} {\bibfnamefont {M.~J.~S.}\ \bibnamefont {Lowe}},\ and\ \bibinfo {author} {\bibfnamefont {P.}~\bibnamefont {Cawley}},\ }\bibfield  {title} {\bibinfo {title} {{Ultrasonic isolation of buried pipes}},\ }\href {https://doi.org/10.1016/j.jsv.2015.10.018} {\bibfield  {journal} {\bibinfo  {journal} {J. Sound Vib.}\ }\textbf {\bibinfo {volume} {363}},\ \bibinfo {pages} {225} (\bibinfo {year} {2016})}\BibitemShut {NoStop}%
\bibitem [{\citenamefont {Leinov}\ \emph {et~al.}(2015)\citenamefont {Leinov}, \citenamefont {Lowe},\ and\ \citenamefont {Cawley}}]{Leinov2015InvestigationSand}%
  \BibitemOpen
  \bibfield  {author} {\bibinfo {author} {\bibfnamefont {E.}~\bibnamefont {Leinov}}, \bibinfo {author} {\bibfnamefont {M.~J.~S.}\ \bibnamefont {Lowe}},\ and\ \bibinfo {author} {\bibfnamefont {P.}~\bibnamefont {Cawley}},\ }\bibfield  {title} {\bibinfo {title} {{Investigation of guided wave propagation and attenuation in pipe buried in sand}},\ }\href {https://doi.org/10.1016/j.jsv.2015.02.036} {\bibfield  {journal} {\bibinfo  {journal} {J. Sound Vib.}\ }\textbf {\bibinfo {volume} {347}},\ \bibinfo {pages} {96} (\bibinfo {year} {2015})}\BibitemShut {NoStop}%
\bibitem [{\citenamefont {Castaings}\ and\ \citenamefont {Lowe}(2008)}]{Castaings2008FiniteMedia}%
  \BibitemOpen
  \bibfield  {author} {\bibinfo {author} {\bibfnamefont {M.}~\bibnamefont {Castaings}}\ and\ \bibinfo {author} {\bibfnamefont {M.~J.~S.}\ \bibnamefont {Lowe}},\ }\bibfield  {title} {\bibinfo {title} {{Finite element model for waves guided along solid systems of arbitrary section coupled to infinite solid media}},\ }\href {https://doi.org/10.1121/1.2821973} {\bibfield  {journal} {\bibinfo  {journal} {J. Acoust. Soc. Am.}\ }\textbf {\bibinfo {volume} {123}},\ \bibinfo {pages} {696} (\bibinfo {year} {2008})}\BibitemShut {NoStop}%
\bibitem [{\citenamefont {Potsika}\ \emph {et~al.}(2014)\citenamefont {Potsika}, \citenamefont {Grivas}, \citenamefont {Protopappas}, \citenamefont {Vavva}, \citenamefont {Raum}, \citenamefont {Rohrbach}, \citenamefont {Polyzos},\ and\ \citenamefont {Fotiadis}}]{Potsika2014ApplicationBones}%
  \BibitemOpen
  \bibfield  {author} {\bibinfo {author} {\bibfnamefont {V.~T.}\ \bibnamefont {Potsika}}, \bibinfo {author} {\bibfnamefont {K.~N.}\ \bibnamefont {Grivas}}, \bibinfo {author} {\bibfnamefont {V.~C.}\ \bibnamefont {Protopappas}}, \bibinfo {author} {\bibfnamefont {M.~G.}\ \bibnamefont {Vavva}}, \bibinfo {author} {\bibfnamefont {K.}~\bibnamefont {Raum}}, \bibinfo {author} {\bibfnamefont {D.}~\bibnamefont {Rohrbach}}, \bibinfo {author} {\bibfnamefont {D.}~\bibnamefont {Polyzos}},\ and\ \bibinfo {author} {\bibfnamefont {D.~I.}\ \bibnamefont {Fotiadis}},\ }\bibfield  {title} {\bibinfo {title} {{Application of an effective medium theory for modeling ultrasound wave propagation in healing long bones}},\ }\href {https://doi.org/10.1016/j.ultras.2013.09.002} {\bibfield  {journal} {\bibinfo  {journal} {Ultrasonics}\ }\textbf {\bibinfo {volume} {54}},\ \bibinfo {pages} {1219} (\bibinfo {year} {2014})}\BibitemShut {NoStop}%
\bibitem [{\citenamefont {Guha}\ \emph {et~al.}(2021)\citenamefont {Guha}, \citenamefont {Aynardi}, \citenamefont {Shokouhi},\ and\ \citenamefont {Lissenden}}]{Guha2021IdentificationModelling}%
  \BibitemOpen
  \bibfield  {author} {\bibinfo {author} {\bibfnamefont {A.}~\bibnamefont {Guha}}, \bibinfo {author} {\bibfnamefont {M.}~\bibnamefont {Aynardi}}, \bibinfo {author} {\bibfnamefont {P.}~\bibnamefont {Shokouhi}},\ and\ \bibinfo {author} {\bibfnamefont {C.~J.}\ \bibnamefont {Lissenden}},\ }\bibfield  {title} {\bibinfo {title} {{Identification of long-range ultrasonic guided wave characteristics in cortical bone by modelling}},\ }\bibfield  {journal} {\bibinfo  {journal} {Ultrasonics}\ }\textbf {\bibinfo {volume} {114}},\ \href {https://doi.org/10.1016/j.ultras.2021.106407} {10.1016/j.ultras.2021.106407} (\bibinfo {year} {2021})\BibitemShut {NoStop}%
\bibitem [{\citenamefont {Lee}\ and\ \citenamefont {Yoon}(2004)}]{Lee2004FeasibilityTibia}%
  \BibitemOpen
  \bibfield  {author} {\bibinfo {author} {\bibfnamefont {K.~I.}\ \bibnamefont {Lee}}\ and\ \bibinfo {author} {\bibfnamefont {S.~W.}\ \bibnamefont {Yoon}},\ }\bibfield  {title} {\bibinfo {title} {{Feasibility of bone assessment with leaky Lamb waves in bone phantoms and a bovine tibia}},\ }\href {https://doi.org/10.1121/1.1707086} {\bibfield  {journal} {\bibinfo  {journal} {J. Acoust. Soc. Am.}\ }\textbf {\bibinfo {volume} {115}},\ \bibinfo {pages} {3210} (\bibinfo {year} {2004})}\BibitemShut {NoStop}%
\bibitem [{\citenamefont {Lamb}(1916)}]{Lamb1916OnPlate}%
  \BibitemOpen
  \bibfield  {author} {\bibinfo {author} {\bibfnamefont {H.}~\bibnamefont {Lamb}},\ }\bibfield  {title} {\bibinfo {title} {{On Waves in an Elastic Plate}},\ }\href {https://doi.org/10.1098/rspa.1917.0008} {\bibfield  {journal} {\bibinfo  {journal} {Proc. R. Soc. A}\ }\textbf {\bibinfo {volume} {81}},\ \bibinfo {pages} {114} (\bibinfo {year} {1916})}\BibitemShut {NoStop}%
\bibitem [{\citenamefont {Kiefer}\ \emph {et~al.}(2019)\citenamefont {Kiefer}, \citenamefont {Ponschab}, \citenamefont {Rupitsch},\ and\ \citenamefont {Mayle}}]{Kiefer2019CalculatingInteraction}%
  \BibitemOpen
  \bibfield  {author} {\bibinfo {author} {\bibfnamefont {D.~A.}\ \bibnamefont {Kiefer}}, \bibinfo {author} {\bibfnamefont {M.}~\bibnamefont {Ponschab}}, \bibinfo {author} {\bibfnamefont {S.~J.}\ \bibnamefont {Rupitsch}},\ and\ \bibinfo {author} {\bibfnamefont {M.}~\bibnamefont {Mayle}},\ }\bibfield  {title} {\bibinfo {title} {{Calculating the full leaky Lamb wave spectrum with exact fluid interaction}},\ }\href {https://doi.org/10.1121/1.5109399} {\bibfield  {journal} {\bibinfo  {journal} {J. Acoust. Soc. Am.}\ }\textbf {\bibinfo {volume} {145}},\ \bibinfo {pages} {3341} (\bibinfo {year} {2019})}\BibitemShut {NoStop}%
\bibitem [{\citenamefont {Pavlakovic}(1998)}]{Pavlakovic1998LeakyNDT}%
  \BibitemOpen
  \bibfield  {author} {\bibinfo {author} {\bibfnamefont {B.~N.}\ \bibnamefont {Pavlakovic}},\ }\emph {\bibinfo {title} {{Leaky Guided Ultrasonic Waves in NDT}}},\ \href@noop {} {Ph.D. thesis},\ \bibinfo  {school} {Imperial College London}, \bibinfo {address} {London} (\bibinfo {year} {1998})\BibitemShut {NoStop}%
\bibitem [{\citenamefont {Georgiades}\ \emph {et~al.}(2022)\citenamefont {Georgiades}, \citenamefont {Lowe},\ and\ \citenamefont {Craster}}]{Georgiades2022LeakyMethods}%
  \BibitemOpen
  \bibfield  {author} {\bibinfo {author} {\bibfnamefont {E.}~\bibnamefont {Georgiades}}, \bibinfo {author} {\bibfnamefont {M.~J.~S.}\ \bibnamefont {Lowe}},\ and\ \bibinfo {author} {\bibfnamefont {R.~V.}\ \bibnamefont {Craster}},\ }\bibfield  {title} {\bibinfo {title} {{Leaky wave characterisation using spectral methods}},\ }\href {https://doi.org/10.1121/10.0013897} {\bibfield  {journal} {\bibinfo  {journal} {J. Acoust. Soc. Am.}\ }\textbf {\bibinfo {volume} {152}},\ \bibinfo {pages} {1487} (\bibinfo {year} {2022})}\BibitemShut {NoStop}%
\bibitem [{\citenamefont {Treyss{\`{e}}de}\ \emph {et~al.}(2014)\citenamefont {Treyss{\`{e}}de}, \citenamefont {Nguyen}, \citenamefont {Bonnet-BenDhia},\ and\ \citenamefont {Hazard}}]{Treyssede2014FiniteWaveguides}%
  \BibitemOpen
  \bibfield  {author} {\bibinfo {author} {\bibfnamefont {F.}~\bibnamefont {Treyss{\`{e}}de}}, \bibinfo {author} {\bibfnamefont {K.~L.}\ \bibnamefont {Nguyen}}, \bibinfo {author} {\bibfnamefont {A.~S.}\ \bibnamefont {Bonnet-BenDhia}},\ and\ \bibinfo {author} {\bibfnamefont {C.}~\bibnamefont {Hazard}},\ }\bibfield  {title} {\bibinfo {title} {{Finite element computation of trapped and leaky elastic waves in open stratified waveguides}},\ }\href {https://doi.org/10.1016/j.wavemoti.2014.05.003} {\bibfield  {journal} {\bibinfo  {journal} {Wave Motion}\ }\textbf {\bibinfo {volume} {51}},\ \bibinfo {pages} {1093} (\bibinfo {year} {2014})}\BibitemShut {NoStop}%
\bibitem [{\citenamefont {Rokhlin}\ \emph {et~al.}(1989)\citenamefont {Rokhlin}, \citenamefont {Chimenti},\ and\ \citenamefont {Nayfeh}}]{Rokhlin1989OnLayer}%
  \BibitemOpen
  \bibfield  {author} {\bibinfo {author} {\bibfnamefont {S.~I.}\ \bibnamefont {Rokhlin}}, \bibinfo {author} {\bibfnamefont {D.~E.}\ \bibnamefont {Chimenti}},\ and\ \bibinfo {author} {\bibfnamefont {A.~H.}\ \bibnamefont {Nayfeh}},\ }\bibfield  {title} {\bibinfo {title} {{On the topology of the complex wave spectrum in a fluid-coupled elastic layer}},\ }\href {https://doi.org/10.1121/1.397490} {\bibfield  {journal} {\bibinfo  {journal} {J. Acoust. Soc. Am.}\ }\textbf {\bibinfo {volume} {85}},\ \bibinfo {pages} {1074} (\bibinfo {year} {1989})}\BibitemShut {NoStop}%
\bibitem [{\citenamefont {Hayashi}\ and\ \citenamefont {Inoue}(2014)}]{Hayashi2014CalculationMethod}%
  \BibitemOpen
  \bibfield  {author} {\bibinfo {author} {\bibfnamefont {T.}~\bibnamefont {Hayashi}}\ and\ \bibinfo {author} {\bibfnamefont {D.}~\bibnamefont {Inoue}},\ }\bibfield  {title} {\bibinfo {title} {{Calculation of leaky Lamb waves with a semi-analytical finite element method}},\ }\href {https://doi.org/10.1016/j.ultras.2014.04.021} {\bibfield  {journal} {\bibinfo  {journal} {Ultrasonics}\ }\textbf {\bibinfo {volume} {54}},\ \bibinfo {pages} {1460} (\bibinfo {year} {2014})}\BibitemShut {NoStop}%
\bibitem [{\citenamefont {Mazzotti}\ \emph {et~al.}(2014)\citenamefont {Mazzotti}, \citenamefont {Bartoli},\ and\ \citenamefont {Marzani}}]{Mazzotti2014UltrasonicValidation}%
  \BibitemOpen
  \bibfield  {author} {\bibinfo {author} {\bibfnamefont {M.}~\bibnamefont {Mazzotti}}, \bibinfo {author} {\bibfnamefont {I.}~\bibnamefont {Bartoli}},\ and\ \bibinfo {author} {\bibfnamefont {A.}~\bibnamefont {Marzani}},\ }\bibfield  {title} {\bibinfo {title} {{Ultrasonic leaky guided waves in fluid-coupled generic waveguides: Hybrid finite-boundary element dispersion analysis and experimental validation}},\ }\bibfield  {journal} {\bibinfo  {journal} {J. Appl. Phys.}\ }\textbf {\bibinfo {volume} {115}},\ \href {https://doi.org/10.1063/1.4870857} {10.1063/1.4870857} (\bibinfo {year} {2014})\BibitemShut {NoStop}%
\bibitem [{\citenamefont {Mazzotti}\ \emph {et~al.}(2013)\citenamefont {Mazzotti}, \citenamefont {Bartoli}, \citenamefont {Marzani},\ and\ \citenamefont {Viola}}]{Mazzotti2013ACross-section}%
  \BibitemOpen
  \bibfield  {author} {\bibinfo {author} {\bibfnamefont {M.}~\bibnamefont {Mazzotti}}, \bibinfo {author} {\bibfnamefont {I.}~\bibnamefont {Bartoli}}, \bibinfo {author} {\bibfnamefont {A.}~\bibnamefont {Marzani}},\ and\ \bibinfo {author} {\bibfnamefont {E.}~\bibnamefont {Viola}},\ }\bibfield  {title} {\bibinfo {title} {{A coupled SAFE-2.5D BEM approach for the dispersion analysis of damped leaky guided waves in embedded waveguides of arbitrary cross-section}},\ }\href {https://doi.org/10.1016/j.ultras.2013.03.003} {\bibfield  {journal} {\bibinfo  {journal} {Ultrasonics}\ }\textbf {\bibinfo {volume} {53}},\ \bibinfo {pages} {1227} (\bibinfo {year} {2013})}\BibitemShut {NoStop}%
\bibitem [{\citenamefont {Zuo}\ \emph {et~al.}(2017)\citenamefont {Zuo}, \citenamefont {Yu},\ and\ \citenamefont {Fan}}]{Zuo2017NumericalPackage}%
  \BibitemOpen
  \bibfield  {author} {\bibinfo {author} {\bibfnamefont {P.}~\bibnamefont {Zuo}}, \bibinfo {author} {\bibfnamefont {X.}~\bibnamefont {Yu}},\ and\ \bibinfo {author} {\bibfnamefont {Z.}~\bibnamefont {Fan}},\ }\bibfield  {title} {\bibinfo {title} {{Numerical modeling of embedded solid waveguides using SAFE-PML approach using a commercially available finite element package}},\ }\href {https://doi.org/10.1016/J.NDTEINT.2017.04.003} {\bibfield  {journal} {\bibinfo  {journal} {NDT E Int.}\ }\textbf {\bibinfo {volume} {90}},\ \bibinfo {pages} {11} (\bibinfo {year} {2017})}\BibitemShut {NoStop}%
\bibitem [{\citenamefont {Boyd}(2000)}]{Boyd2000ChebyshevMethods}%
  \BibitemOpen
  \bibfield  {author} {\bibinfo {author} {\bibfnamefont {J.~P.}\ \bibnamefont {Boyd}},\ }\href@noop {} {\emph {\bibinfo {title} {{Chebyshev and Fourier Spectral Methods}}}},\ \bibinfo {edition} {second edition}\ ed.\ (\bibinfo  {publisher} {Dover Publications},\ \bibinfo {address} {Mineola, New York},\ \bibinfo {year} {2000})\ pp.\ \bibinfo {pages} {1--594}\BibitemShut {NoStop}%
\bibitem [{\citenamefont {Trefethen}(2008)}]{Trefethen2008SpectralMatlab}%
  \BibitemOpen
  \bibfield  {author} {\bibinfo {author} {\bibfnamefont {L.~N.}\ \bibnamefont {Trefethen}},\ }\href {https://doi.org/10.1137/1.9780898719598} {\emph {\bibinfo {title} {{Spectral Methods in Matlab}}}}\ (\bibinfo  {publisher} {SIAM},\ \bibinfo {address} {Philadelphia},\ \bibinfo {year} {2008})\ pp.\ \bibinfo {pages} {1--160}\BibitemShut {NoStop}%
\bibitem [{\citenamefont {Ryskamp}\ and\ \citenamefont {Hoefer}(2023)}]{Ryskamp2023ComputationMaterials}%
  \BibitemOpen
  \bibfield  {author} {\bibinfo {author} {\bibfnamefont {S.~J.}\ \bibnamefont {Ryskamp}}\ and\ \bibinfo {author} {\bibfnamefont {M.~A.}\ \bibnamefont {Hoefer}},\ }\bibfield  {title} {\bibinfo {title} {{Computation of High Frequency Magnetoelastic Waves in Layered Materials}},\ }\bibfield  {journal} {\bibinfo  {journal} {Phys. Rev. B}\ }\textbf {\bibinfo {volume} {107}},\ \href {https://doi.org/10.1103/PhysRevB.107.094431} {10.1103/PhysRevB.107.094431} (\bibinfo {year} {2023})\BibitemShut {NoStop}%
\bibitem [{\citenamefont {Adamou}\ and\ \citenamefont {Craster}(2004)}]{Adamou2004SpectralMedia}%
  \BibitemOpen
  \bibfield  {author} {\bibinfo {author} {\bibfnamefont {A.~T.~I.}\ \bibnamefont {Adamou}}\ and\ \bibinfo {author} {\bibfnamefont {R.~V.}\ \bibnamefont {Craster}},\ }\bibfield  {title} {\bibinfo {title} {{Spectral methods for modelling guided waves in elastic media}},\ }\href {https://doi.org/10.1121/1.1777871} {\bibfield  {journal} {\bibinfo  {journal} {J. Acoust. Soc. Am.}\ }\textbf {\bibinfo {volume} {116}},\ \bibinfo {pages} {1524} (\bibinfo {year} {2004})}\BibitemShut {NoStop}%
\bibitem [{\citenamefont {Quintanilla}\ \emph {et~al.}(2016)\citenamefont {Quintanilla}, \citenamefont {Lowe},\ and\ \citenamefont {Craster}}]{Quintanilla2016FullMedia}%
  \BibitemOpen
  \bibfield  {author} {\bibinfo {author} {\bibfnamefont {F.~H.}\ \bibnamefont {Quintanilla}}, \bibinfo {author} {\bibfnamefont {M.~J.~S.}\ \bibnamefont {Lowe}},\ and\ \bibinfo {author} {\bibfnamefont {R.~V.}\ \bibnamefont {Craster}},\ }\bibfield  {title} {\bibinfo {title} {{Full 3D dispersion curve solutions for guided waves in generally anisotropic media}},\ }\href {https://doi.org/10.1016/j.jsv.2015.10.017} {\bibfield  {journal} {\bibinfo  {journal} {J. Sound Vib.}\ }\textbf {\bibinfo {volume} {363}},\ \bibinfo {pages} {545} (\bibinfo {year} {2016})}\BibitemShut {NoStop}%
\bibitem [{\citenamefont {Quintanilla}\ \emph {et~al.}(2015{\natexlab{a}})\citenamefont {Quintanilla}, \citenamefont {Fan}, \citenamefont {Lowe},\ and\ \citenamefont {Craster}}]{Quintanilla2015GuidedMedia}%
  \BibitemOpen
  \bibfield  {author} {\bibinfo {author} {\bibfnamefont {F.~H.}\ \bibnamefont {Quintanilla}}, \bibinfo {author} {\bibfnamefont {Z.}~\bibnamefont {Fan}}, \bibinfo {author} {\bibfnamefont {M.~J.~S.}\ \bibnamefont {Lowe}},\ and\ \bibinfo {author} {\bibfnamefont {R.~V.}\ \bibnamefont {Craster}},\ }\bibfield  {title} {\bibinfo {title} {{Guided waves' dispersion curves in anisotropic viscoelastic single- and multi-layered media}},\ }\bibfield  {journal} {\bibinfo  {journal} {Proc. R. Soc. A}\ }\textbf {\bibinfo {volume} {471}},\ \href {https://doi.org/10.1098/rspa.2015.0268} {10.1098/rspa.2015.0268} (\bibinfo {year} {2015}{\natexlab{a}})\BibitemShut {NoStop}%
\bibitem [{\citenamefont {Quintanilla}\ \emph {et~al.}(2015{\natexlab{b}})\citenamefont {Quintanilla}, \citenamefont {Lowe},\ and\ \citenamefont {Craster}}]{Quintanilla2015ModelingMethod}%
  \BibitemOpen
  \bibfield  {author} {\bibinfo {author} {\bibfnamefont {F.~H.}\ \bibnamefont {Quintanilla}}, \bibinfo {author} {\bibfnamefont {M.~J.~S.}\ \bibnamefont {Lowe}},\ and\ \bibinfo {author} {\bibfnamefont {R.~V.}\ \bibnamefont {Craster}},\ }\bibfield  {title} {\bibinfo {title} {{Modeling guided elastic waves in generally anisotropic media using a spectral collocation method}},\ }\href {https://doi.org/10.1121/1.4913777} {\bibfield  {journal} {\bibinfo  {journal} {J. Acoust. Soc. Am.}\ }\textbf {\bibinfo {volume} {137}},\ \bibinfo {pages} {1180} (\bibinfo {year} {2015}{\natexlab{b}})}\BibitemShut {NoStop}%
\bibitem [{\citenamefont {Quintanilla}\ \emph {et~al.}(2017)\citenamefont {Quintanilla}, \citenamefont {Lowe},\ and\ \citenamefont {Craster}}]{Quintanilla2017TheWaveguides}%
  \BibitemOpen
  \bibfield  {author} {\bibinfo {author} {\bibfnamefont {F.~H.}\ \bibnamefont {Quintanilla}}, \bibinfo {author} {\bibfnamefont {M.~J.~S.}\ \bibnamefont {Lowe}},\ and\ \bibinfo {author} {\bibfnamefont {R.~V.}\ \bibnamefont {Craster}},\ }\bibfield  {title} {\bibinfo {title} {{The symmetry and coupling properties of solutions in general anisotropic multilayer waveguides}},\ }\href {https://doi.org/10.1121/1.4973543} {\bibfield  {journal} {\bibinfo  {journal} {J. Acoust. Soc. Am.}\ }\textbf {\bibinfo {volume} {141}},\ \bibinfo {pages} {406} (\bibinfo {year} {2017})}\BibitemShut {NoStop}%
\bibitem [{\citenamefont {Tang}\ \emph {et~al.}(2022)\citenamefont {Tang}, \citenamefont {Yin}, \citenamefont {Wang},\ and\ \citenamefont {Zhu}}]{Tang2022StudyFluids}%
  \BibitemOpen
  \bibfield  {author} {\bibinfo {author} {\bibfnamefont {S.}~\bibnamefont {Tang}}, \bibinfo {author} {\bibfnamefont {J.}~\bibnamefont {Yin}}, \bibinfo {author} {\bibfnamefont {C.}~\bibnamefont {Wang}},\ and\ \bibinfo {author} {\bibfnamefont {G.}~\bibnamefont {Zhu}},\ }\bibfield  {title} {\bibinfo {title} {{Study on leaky Lamb waves in functionally graded composites loaded by asymmetric fluids}},\ }\bibfield  {journal} {\bibinfo  {journal} {Waves Random Complex Media}\ }\href {https://doi.org/10.1080/17455030.2022.2102270} {10.1080/17455030.2022.2102270} (\bibinfo {year} {2022})\BibitemShut {NoStop}%
\bibitem [{\citenamefont {Pavlakovic}\ \emph {et~al.}(1997)\citenamefont {Pavlakovic}, \citenamefont {Lowe}, \citenamefont {Alleyne},\ and\ \citenamefont {Cawley}}]{Pavlakovic1997Disperse:Curves}%
  \BibitemOpen
  \bibfield  {author} {\bibinfo {author} {\bibfnamefont {B.}~\bibnamefont {Pavlakovic}}, \bibinfo {author} {\bibfnamefont {M.~J.~S.}\ \bibnamefont {Lowe}}, \bibinfo {author} {\bibfnamefont {O.}~\bibnamefont {Alleyne}},\ and\ \bibinfo {author} {\bibfnamefont {P.}~\bibnamefont {Cawley}},\ }\bibfield  {title} {\bibinfo {title} {{Disperse: A general purpose program for creating dispersion curves}},\ }\href {https://doi.org/10.1007/978-1-4615-5947-4{\_}24} {\bibfield  {journal} {\bibinfo  {journal} {Review of Progress in QNDE}\ }\textbf {\bibinfo {volume} {16}},\ \bibinfo {pages} {185} (\bibinfo {year} {1997})}\BibitemShut {NoStop}%
\bibitem [{\citenamefont {{Bathe K. J.}}(1996)}]{BatheK.J.1996FiniteProcedures}%
  \BibitemOpen
  \bibfield  {author} {\bibinfo {author} {\bibnamefont {{Bathe K. J.}}},\ }\href@noop {} {\emph {\bibinfo {title} {{Finite Element Procedures}}}},\ \bibinfo {edition} {1st}\ ed.\ (\bibinfo  {publisher} {Prentice Hal},\ \bibinfo {year} {1996})\BibitemShut {NoStop}%
\bibitem [{\citenamefont {Auld}(1973)}]{Auld1973AcousticSolidsb}%
  \BibitemOpen
  \bibfield  {author} {\bibinfo {author} {\bibfnamefont {B.~A.}\ \bibnamefont {Auld}},\ }\href@noop {} {\emph {\bibinfo {title} {{Acoustic fields and waves in solids}}}},\ Vol.~\bibinfo {volume} {2}\ (\bibinfo {year} {1973})\BibitemShut {NoStop}%
\bibitem [{\citenamefont {Nayfeh}\ and\ \citenamefont {Nagy}(1997)}]{Nayfeh1997ExcessLoading}%
  \BibitemOpen
  \bibfield  {author} {\bibinfo {author} {\bibfnamefont {A.~H.}\ \bibnamefont {Nayfeh}}\ and\ \bibinfo {author} {\bibfnamefont {P.~B.}\ \bibnamefont {Nagy}},\ }\bibfield  {title} {\bibinfo {title} {{Excess attenuation of leaky Lamb waves due to viscous fluid loading}},\ }\href {https://doi.org/10.1121/1.418506} {\bibfield  {journal} {\bibinfo  {journal} {J. Acoust. Soc. Am.}\ }\textbf {\bibinfo {volume} {101}},\ \bibinfo {pages} {2649} (\bibinfo {year} {1997})}\BibitemShut {NoStop}%
\bibitem [{\citenamefont {Achenbach}(1973)}]{Achenbach1973WaveSolids}%
  \BibitemOpen
  \bibfield  {author} {\bibinfo {author} {\bibfnamefont {J.~D.}\ \bibnamefont {Achenbach}},\ }\href@noop {} {\emph {\bibinfo {title} {{Wave propagation in elastic solids}}}}\ (\bibinfo  {publisher} {North-Holland Pub. Co},\ \bibinfo {year} {1973})\ p.\ \bibinfo {pages} {425}\BibitemShut {NoStop}%
\bibitem [{\citenamefont {Rose}(2014)}]{Rose2014UltrasonicMedia}%
  \BibitemOpen
  \bibfield  {author} {\bibinfo {author} {\bibfnamefont {J.~L.}\ \bibnamefont {Rose}},\ }\href {https://doi.org/10.1017/CBO9781107273610} {\emph {\bibinfo {title} {Ultrasonic Guided Waves in Solid Media}}}\ (\bibinfo  {publisher} {Cambridge University Press},\ \bibinfo {year} {2014})\ pp.\ \bibinfo {pages} {1--512}\BibitemShut {NoStop}%
\bibitem [{\citenamefont {Gallezot}\ \emph {et~al.}(2017)\citenamefont {Gallezot}, \citenamefont {Treyss{\`{e}}de},\ and\ \citenamefont {Laguerre}}]{Gallezot2017ContributionLayers}%
  \BibitemOpen
  \bibfield  {author} {\bibinfo {author} {\bibfnamefont {M.}~\bibnamefont {Gallezot}}, \bibinfo {author} {\bibfnamefont {F.}~\bibnamefont {Treyss{\`{e}}de}},\ and\ \bibinfo {author} {\bibfnamefont {L.}~\bibnamefont {Laguerre}},\ }\bibfield  {title} {\bibinfo {title} {{Contribution of leaky modes in the modal analysis of unbounded problems with perfectly matched layers}},\ }\href {https://doi.org/10.1121/1.4973313} {\bibfield  {journal} {\bibinfo  {journal} {J. Acoust. Soc. Am.}\ }\textbf {\bibinfo {volume} {141}},\ \bibinfo {pages} {EL16} (\bibinfo {year} {2017})}\BibitemShut {NoStop}%
\bibitem [{\citenamefont {Kiefer}(2022)}]{Kiefer2022ElastodynamicMetering}%
  \BibitemOpen
  \bibfield  {author} {\bibinfo {author} {\bibfnamefont {D.~A.}\ \bibnamefont {Kiefer}},\ }\emph {\bibinfo {title} {{Elastodynamic quasi-guided waves for transit-time ultrasonic flow metering}}},\ \href {https://doi.org/10.25593/978-3-96147-550-6} {Ph.D. thesis},\ \bibinfo  {school} {Universitaet Erlangen-Nuernberg}, \bibinfo {address} {Erlangen} (\bibinfo {year} {2022})\BibitemShut {NoStop}%
\bibitem [{\citenamefont {Mozhaev}\ and\ \citenamefont {Weihnacht}(2002)}]{Mozhaev2002SubsonicInterfaces}%
  \BibitemOpen
  \bibfield  {author} {\bibinfo {author} {\bibfnamefont {V.~G.}\ \bibnamefont {Mozhaev}}\ and\ \bibinfo {author} {\bibfnamefont {M.}~\bibnamefont {Weihnacht}},\ }\bibfield  {title} {\bibinfo {title} {{Subsonic leaky Rayleigh waves at liquid-solid interfaces}},\ }\href {https://doi.org/10.1016/s0041-624x(02)00233-0} {\bibfield  {journal} {\bibinfo  {journal} {Ultrasonics}\ }\textbf {\bibinfo {volume} {40}},\ \bibinfo {pages} {927} (\bibinfo {year} {2002})}\BibitemShut {NoStop}%
\bibitem [{\citenamefont {Viggen}\ and\ \citenamefont {Arnestad}(2023)}]{Viggen2023ModellingFluids}%
  \BibitemOpen
  \bibfield  {author} {\bibinfo {author} {\bibfnamefont {E.~M.}\ \bibnamefont {Viggen}}\ and\ \bibinfo {author} {\bibfnamefont {H.~K.}\ \bibnamefont {Arnestad}},\ }\bibfield  {title} {\bibinfo {title} {{Modelling acoustic radiation from vibrating surfaces around coincidence: Radiation into fluids}},\ }\href {https://doi.org/10.1016/j.jsv.2023.117787} {\bibfield  {journal} {\bibinfo  {journal} {J. Sound Vib.}\ }\textbf {\bibinfo {volume} {560}},\ \bibinfo {pages} {117787} (\bibinfo {year} {2023})}\BibitemShut {NoStop}%
\bibitem [{\citenamefont {Fornberg}(1996)}]{Fornberg1996AMethods}%
  \BibitemOpen
  \bibfield  {author} {\bibinfo {author} {\bibfnamefont {B.}~\bibnamefont {Fornberg}},\ }\href {https://doi.org/10.1017/CBO9780511626357} {\emph {\bibinfo {title} {A Practical Guide to Pseudospectral Methods}}}\ (\bibinfo  {publisher} {Cambridge University Press},\ \bibinfo {address} {Cambridge},\ \bibinfo {year} {1996})\ pp.\ \bibinfo {pages} {1--231}\BibitemShut {NoStop}%
\bibitem [{\citenamefont {Weideman}\ and\ \citenamefont {Reddy}(2000)}]{Weideman2000ASuite}%
  \BibitemOpen
  \bibfield  {author} {\bibinfo {author} {\bibfnamefont {J.~A.~C.}\ \bibnamefont {Weideman}}\ and\ \bibinfo {author} {\bibfnamefont {S.~C.}\ \bibnamefont {Reddy}},\ }\bibfield  {title} {\bibinfo {title} {{A MATLAB Differentiation Matrix Suite}},\ }\href {https://doi.org/10.1145/365723.365727} {\bibfield  {journal} {\bibinfo  {journal} {ACM Trans. Math. Softw.}\ }\textbf {\bibinfo {volume} {26}},\ \bibinfo {pages} {465} (\bibinfo {year} {2000})}\BibitemShut {NoStop}%
\bibitem [{\citenamefont {Mackey}\ \emph {et~al.}(2006)\citenamefont {Mackey}, \citenamefont {Mackey}, \citenamefont {Mehl},\ and\ \citenamefont {Mehrmann}}]{Mackey2006StructuredLinearizations}%
  \BibitemOpen
  \bibfield  {author} {\bibinfo {author} {\bibfnamefont {D.~S.}\ \bibnamefont {Mackey}}, \bibinfo {author} {\bibfnamefont {N.}~\bibnamefont {Mackey}}, \bibinfo {author} {\bibfnamefont {C.}~\bibnamefont {Mehl}},\ and\ \bibinfo {author} {\bibfnamefont {V.}~\bibnamefont {Mehrmann}},\ }\bibfield  {title} {\bibinfo {title} {{Structured polynomial eigenvalue problems: Good vibrations from good linearizations}},\ }\href {https://doi.org/10.1137/050628362} {\bibfield  {journal} {\bibinfo  {journal} {SIAM J. Matrix Anal. Appl.}\ }\textbf {\bibinfo {volume} {28}},\ \bibinfo {pages} {1029} (\bibinfo {year} {2006})}\BibitemShut {NoStop}%
\bibitem [{\citenamefont {Huthwaite}(2014)}]{Huthwaite2014AcceleratedGPU}%
  \BibitemOpen
  \bibfield  {author} {\bibinfo {author} {\bibfnamefont {P.}~\bibnamefont {Huthwaite}},\ }\bibfield  {title} {\bibinfo {title} {{Accelerated finite element elastodynamic simulations using the GPU}},\ }\href {https://doi.org/10.1016/j.jcp.2013.10.017} {\bibfield  {journal} {\bibinfo  {journal} {J. Comput. Phys.}\ }\textbf {\bibinfo {volume} {257}},\ \bibinfo {pages} {687} (\bibinfo {year} {2014})}\BibitemShut {NoStop}%
\bibitem [{\citenamefont {Sarris}\ \emph {et~al.}(2021)\citenamefont {Sarris}, \citenamefont {Haslinger}, \citenamefont {Huthwaite}, \citenamefont {Nagy},\ and\ \citenamefont {Lowe}}]{Sarris2021AttenuationRoughness}%
  \BibitemOpen
  \bibfield  {author} {\bibinfo {author} {\bibfnamefont {G.}~\bibnamefont {Sarris}}, \bibinfo {author} {\bibfnamefont {S.~G.}\ \bibnamefont {Haslinger}}, \bibinfo {author} {\bibfnamefont {P.}~\bibnamefont {Huthwaite}}, \bibinfo {author} {\bibfnamefont {P.~B.}\ \bibnamefont {Nagy}},\ and\ \bibinfo {author} {\bibfnamefont {M.~J.~S.}\ \bibnamefont {Lowe}},\ }\bibfield  {title} {\bibinfo {title} {{Attenuation of Rayleigh waves due to surface roughness}},\ }\href {https://doi.org/10.1121/10.0005271} {\bibfield  {journal} {\bibinfo  {journal} {J. Acoust. Soc. Am.}\ }\textbf {\bibinfo {volume} {149}},\ \bibinfo {pages} {4298} (\bibinfo {year} {2021})}\BibitemShut {NoStop}%
\bibitem [{\citenamefont {Sarris}\ \emph {et~al.}(2023)\citenamefont {Sarris}, \citenamefont {Haslinger}, \citenamefont {Huthwaite},\ and\ \citenamefont {Lowe}}]{Sarris2023UltrasonicDamage}%
  \BibitemOpen
  \bibfield  {author} {\bibinfo {author} {\bibfnamefont {G.}~\bibnamefont {Sarris}}, \bibinfo {author} {\bibfnamefont {S.~G.}\ \bibnamefont {Haslinger}}, \bibinfo {author} {\bibfnamefont {P.}~\bibnamefont {Huthwaite}},\ and\ \bibinfo {author} {\bibfnamefont {M.~J.}\ \bibnamefont {Lowe}},\ }\bibfield  {title} {\bibinfo {title} {{Ultrasonic methods for the detection of near surface fatigue damage}},\ }\bibfield  {journal} {\bibinfo  {journal} {NDT E Int.}\ }\textbf {\bibinfo {volume} {135}},\ \href {https://doi.org/10.1016/j.ndteint.2023.102790} {10.1016/j.ndteint.2023.102790} (\bibinfo {year} {2023})\BibitemShut {NoStop}%
\bibitem [{\citenamefont {Haslinger}\ \emph {et~al.}(2019)\citenamefont {Haslinger}, \citenamefont {Lowe}, \citenamefont {Huthwaite}, \citenamefont {Craster},\ and\ \citenamefont {Shi}}]{Haslinger2019AppraisingDefects}%
  \BibitemOpen
  \bibfield  {author} {\bibinfo {author} {\bibfnamefont {S.~G.}\ \bibnamefont {Haslinger}}, \bibinfo {author} {\bibfnamefont {M.~J.}\ \bibnamefont {Lowe}}, \bibinfo {author} {\bibfnamefont {P.}~\bibnamefont {Huthwaite}}, \bibinfo {author} {\bibfnamefont {R.~V.}\ \bibnamefont {Craster}},\ and\ \bibinfo {author} {\bibfnamefont {F.}~\bibnamefont {Shi}},\ }\bibfield  {title} {\bibinfo {title} {{Appraising Kirchhoff approximation theory for the scattering of elastic shear waves by randomly rough defects}},\ }\bibfield  {journal} {\bibinfo  {journal} {J. Sound Vib.}\ }\textbf {\bibinfo {volume} {460}},\ \href {https://doi.org/10.1016/j.jsv.2019.114872} {10.1016/j.jsv.2019.114872} (\bibinfo {year} {2019})\BibitemShut {NoStop}%
\bibitem [{\citenamefont {Huang}\ \emph {et~al.}(2020)\citenamefont {Huang}, \citenamefont {Sha}, \citenamefont {Huthwaite}, \citenamefont {Rokhlin},\ and\ \citenamefont {Lowe}}]{Huang2020MaximizingPolycrystals}%
  \BibitemOpen
  \bibfield  {author} {\bibinfo {author} {\bibfnamefont {M.}~\bibnamefont {Huang}}, \bibinfo {author} {\bibfnamefont {G.}~\bibnamefont {Sha}}, \bibinfo {author} {\bibfnamefont {P.}~\bibnamefont {Huthwaite}}, \bibinfo {author} {\bibfnamefont {S.~I.}\ \bibnamefont {Rokhlin}},\ and\ \bibinfo {author} {\bibfnamefont {M.~J.~S.}\ \bibnamefont {Lowe}},\ }\bibfield  {title} {\bibinfo {title} {{Maximizing the accuracy of finite element simulation of elastic wave propagation in polycrystals}},\ }\href {https://doi.org/10.1121/10.0002102} {\bibfield  {journal} {\bibinfo  {journal} {J. Acoust. Soc. Am.}\ }\textbf {\bibinfo {volume} {148}},\ \bibinfo {pages} {1890} (\bibinfo {year} {2020})}\BibitemShut {NoStop}%
\bibitem [{\citenamefont {Drozdz}(2008)}]{Drozdz2008EfficientMedia}%
  \BibitemOpen
  \bibfield  {author} {\bibinfo {author} {\bibfnamefont {M.~B.}\ \bibnamefont {Drozdz}},\ }\emph {\bibinfo {title} {{Efficient Finite Element Modelling of Ultrasonic Waves in Elastic Media}}},\ \href@noop {} {Ph.D. thesis},\ \bibinfo  {school} {Imperial College London}, \bibinfo {address} {London} (\bibinfo {year} {2008})\BibitemShut {NoStop}%
\bibitem [{\citenamefont {Rajagopal}\ \emph {et~al.}(2012)\citenamefont {Rajagopal}, \citenamefont {Drozdz}, \citenamefont {Skelton}, \citenamefont {Lowe},\ and\ \citenamefont {Craster}}]{Rajagopal2012OnPackages}%
  \BibitemOpen
  \bibfield  {author} {\bibinfo {author} {\bibfnamefont {P.}~\bibnamefont {Rajagopal}}, \bibinfo {author} {\bibfnamefont {M.~B.}\ \bibnamefont {Drozdz}}, \bibinfo {author} {\bibfnamefont {E.~A.}\ \bibnamefont {Skelton}}, \bibinfo {author} {\bibfnamefont {M.~J.}\ \bibnamefont {Lowe}},\ and\ \bibinfo {author} {\bibfnamefont {R.~V.}\ \bibnamefont {Craster}},\ }\bibfield  {title} {\bibinfo {title} {{On the use of absorbing layers to simulate the propagation of elastic waves in unbounded isotropic media using commercially available Finite Element packages}},\ }\href {https://doi.org/10.1016/j.ndteint.2012.04.001} {\bibfield  {journal} {\bibinfo  {journal} {NDT E Int.}\ }\textbf {\bibinfo {volume} {51}},\ \bibinfo {pages} {30} (\bibinfo {year} {2012})}\BibitemShut {NoStop}%
\bibitem [{\citenamefont {Gallezot}\ \emph {et~al.}(2018)\citenamefont {Gallezot}, \citenamefont {Treyss{\`{e}}de},\ and\ \citenamefont {Laguerre}}]{Gallezot2018AWaveguides}%
  \BibitemOpen
  \bibfield  {author} {\bibinfo {author} {\bibfnamefont {M.}~\bibnamefont {Gallezot}}, \bibinfo {author} {\bibfnamefont {F.}~\bibnamefont {Treyss{\`{e}}de}},\ and\ \bibinfo {author} {\bibfnamefont {L.}~\bibnamefont {Laguerre}},\ }\bibfield  {title} {\bibinfo {title} {{A modal approach based on perfectly matched layers for the forced response of elastic open waveguides}},\ }\href {https://doi.org/10.1016/j.jcp.2017.12.017} {\bibfield  {journal} {\bibinfo  {journal} {J. Comput. Phys.}\ }\textbf {\bibinfo {volume} {356}},\ \bibinfo {pages} {391} (\bibinfo {year} {2018})}\BibitemShut {NoStop}%
\bibitem [{\citenamefont {Sammut}\ and\ \citenamefont {Snyder}(1976)}]{Sammut1976LeakyExcitation}%
  \BibitemOpen
  \bibfield  {author} {\bibinfo {author} {\bibfnamefont {R.}~\bibnamefont {Sammut}}\ and\ \bibinfo {author} {\bibfnamefont {A.~W.}\ \bibnamefont {Snyder}},\ }\bibfield  {title} {\bibinfo {title} {{Leaky modes on a dielectric waveguide: orthogonality and excitation}},\ }\href {https://doi.org/10.1364/AO.15.001040} {\bibfield  {journal} {\bibinfo  {journal} {Appl. Opt.}\ }\textbf {\bibinfo {volume} {15}},\ \bibinfo {pages} {1040} (\bibinfo {year} {1976})}\BibitemShut {NoStop}%
\bibitem [{\citenamefont {Snyder}\ and\ \citenamefont {Love}(1983)}]{Snyder1983OpticalTheory}%
  \BibitemOpen
  \bibfield  {author} {\bibinfo {author} {\bibfnamefont {A.~W.}\ \bibnamefont {Snyder}}\ and\ \bibinfo {author} {\bibfnamefont {J.}~\bibnamefont {Love}},\ }\href {https://doi.org/10.1007/978-1-4613-2813-1} {\emph {\bibinfo {title} {{Optical Waveguide Theory}}}}\ (\bibinfo  {publisher} {Chapman and Hall},\ \bibinfo {address} {London},\ \bibinfo {year} {1983})\BibitemShut {NoStop}%
\bibitem [{\citenamefont {Alleyne}\ and\ \citenamefont {Cawley}(1991)}]{Alleyne1991ASignals}%
  \BibitemOpen
  \bibfield  {author} {\bibinfo {author} {\bibfnamefont {D.}~\bibnamefont {Alleyne}}\ and\ \bibinfo {author} {\bibfnamefont {P.}~\bibnamefont {Cawley}},\ }\bibfield  {title} {\bibinfo {title} {{A two-dimensional Fourier transform method for the measurement of propagating multimode signals}},\ }\href {https://doi.org/10.1121/1.400530} {\bibfield  {journal} {\bibinfo  {journal} {J. Acoust. Soc. Am.}\ }\textbf {\bibinfo {volume} {89}},\ \bibinfo {pages} {1159} (\bibinfo {year} {1991})}\BibitemShut {NoStop}%
\bibitem [{\citenamefont {Lowe}\ and\ \citenamefont {Cawley}(1994)}]{Lowe1994TheJoints}%
  \BibitemOpen
  \bibfield  {author} {\bibinfo {author} {\bibfnamefont {M.~J.~S.}\ \bibnamefont {Lowe}}\ and\ \bibinfo {author} {\bibfnamefont {P.}~\bibnamefont {Cawley}},\ }\bibfield  {title} {\bibinfo {title} {{The Applicability of Plate Wave Techniques for the Inspection of Adhesive and Diffusion Bonded Joints}},\ }\href {https://doi.org/10.1007/BF00742584} {\bibfield  {journal} {\bibinfo  {journal} {J. Nondestruct. Eval.}\ }\textbf {\bibinfo {volume} {13}},\ \bibinfo {pages} {185} (\bibinfo {year} {1994})}\BibitemShut {NoStop}%
\bibitem [{\citenamefont {Berenger}(1994)}]{Berenger1994AWaves}%
  \BibitemOpen
  \bibfield  {author} {\bibinfo {author} {\bibfnamefont {J.-P.}\ \bibnamefont {Berenger}},\ }\bibfield  {title} {\bibinfo {title} {{A Perfectly Matched Layer for the Absorption of Electromagnetic Waves}},\ }\href {https://doi.org/10.1006/jcph.1994.1159} {\bibfield  {journal} {\bibinfo  {journal} {J. Comput. Phys}\ }\textbf {\bibinfo {volume} {114}},\ \bibinfo {pages} {185} (\bibinfo {year} {1994})}\BibitemShut {NoStop}%
\bibitem [{\citenamefont {Chew}\ and\ \citenamefont {Weedon}(1994)}]{Chew1994ACoordinates}%
  \BibitemOpen
  \bibfield  {author} {\bibinfo {author} {\bibfnamefont {W.~C.}\ \bibnamefont {Chew}}\ and\ \bibinfo {author} {\bibfnamefont {W.~H.}\ \bibnamefont {Weedon}},\ }\bibfield  {title} {\bibinfo {title} {{A 3D Perfectly matched medium from modified Maxwell's equations with stretched coordinates}},\ }\href {https://doi.org/10.1002/mop.4650071304} {\bibfield  {journal} {\bibinfo  {journal} {Microw. Opt. Technol. Lett.}\ }\textbf {\bibinfo {volume} {7}},\ \bibinfo {pages} {599} (\bibinfo {year} {1994})}\BibitemShut {NoStop}%
\bibitem [{\citenamefont {Chew}\ and\ \citenamefont {Liu}(1996)}]{Chew1996PerfectlyCondition}%
  \BibitemOpen
  \bibfield  {author} {\bibinfo {author} {\bibfnamefont {W.~C.}\ \bibnamefont {Chew}}\ and\ \bibinfo {author} {\bibfnamefont {Q.~H.}\ \bibnamefont {Liu}},\ }\bibfield  {title} {\bibinfo {title} {{Perfectly Matched Layers for Elastodynamics: A New Absorbing Boundary Condition}},\ }\href {https://doi.org/10.1142/S0218396X96000118} {\bibfield  {journal} {\bibinfo  {journal} {J. Comput. Acoust.}\ }\textbf {\bibinfo {volume} {4}} (\bibinfo {year} {1996})}\BibitemShut {NoStop}%
\bibitem [{\citenamefont {Basu}\ and\ \citenamefont {Chopra}(2003)}]{Basu2003PerfectlyImplementation}%
  \BibitemOpen
  \bibfield  {author} {\bibinfo {author} {\bibfnamefont {U.}~\bibnamefont {Basu}}\ and\ \bibinfo {author} {\bibfnamefont {A.~K.}\ \bibnamefont {Chopra}},\ }\bibfield  {title} {\bibinfo {title} {{Perfectly matched layers for time-harmonic elastodynamics of unbounded domains: Theory and finite-element implementation}},\ }\href {https://doi.org/10.1016/S0045-7825(02)00642-4} {\bibfield  {journal} {\bibinfo  {journal} {Comput. Methods Appl. Mech. Eng.}\ }\textbf {\bibinfo {volume} {192}},\ \bibinfo {pages} {1337} (\bibinfo {year} {2003})}\BibitemShut {NoStop}%
\bibitem [{\citenamefont {Kim}\ and\ \citenamefont {Pasciak}(2009)}]{Kim2009TheLayer}%
  \BibitemOpen
  \bibfield  {author} {\bibinfo {author} {\bibfnamefont {S.}~\bibnamefont {Kim}}\ and\ \bibinfo {author} {\bibfnamefont {J.~E.}\ \bibnamefont {Pasciak}},\ }\bibfield  {title} {\bibinfo {title} {{The Computation of Resonances in Open Systems Using a Perfectly Matched Layer}},\ }\href {https://doi.org/10.1090/S0025-5718-09-02227-3} {\bibfield  {journal} {\bibinfo  {journal} {Math. Comp.}\ }\textbf {\bibinfo {volume} {78}},\ \bibinfo {pages} {1375} (\bibinfo {year} {2009})}\BibitemShut {NoStop}%
\bibitem [{\citenamefont {Collino}\ and\ \citenamefont {Tsogka}(2001)}]{Collino2001ApplicationMedia}%
  \BibitemOpen
  \bibfield  {author} {\bibinfo {author} {\bibfnamefont {F.}~\bibnamefont {Collino}}\ and\ \bibinfo {author} {\bibfnamefont {C.}~\bibnamefont {Tsogka}},\ }\bibfield  {title} {\bibinfo {title} {{Application of the perfectly matched absorbing layer model to the linear elastodynamic problem in anisotropic heterogeneous media}},\ }\href {https://doi.org/10.1190/1.1444908} {\bibfield  {journal} {\bibinfo  {journal} {Geophysics}\ }\textbf {\bibinfo {volume} {66}},\ \bibinfo {pages} {294} (\bibinfo {year} {2001})}\BibitemShut {NoStop}%
\bibitem [{\citenamefont {Harari}\ and\ \citenamefont {Albocher}(2006)}]{Harari2006StudiesWaves}%
  \BibitemOpen
  \bibfield  {author} {\bibinfo {author} {\bibfnamefont {I.}~\bibnamefont {Harari}}\ and\ \bibinfo {author} {\bibfnamefont {U.}~\bibnamefont {Albocher}},\ }\bibfield  {title} {\bibinfo {title} {{Studies of FE/PML for exterior problems of time-harmonic elastic waves}},\ }\href {https://doi.org/10.1016/j.cma.2005.01.024} {\bibfield  {journal} {\bibinfo  {journal} {Comput. Methods Appl. Mech. Eng.}\ }\textbf {\bibinfo {volume} {195}},\ \bibinfo {pages} {3854} (\bibinfo {year} {2006})}\BibitemShut {NoStop}%
\bibitem [{\citenamefont {Pled}\ and\ \citenamefont {Desceliers}(2022)}]{Pled2022ReviewDomains}%
  \BibitemOpen
  \bibfield  {author} {\bibinfo {author} {\bibfnamefont {F.}~\bibnamefont {Pled}}\ and\ \bibinfo {author} {\bibfnamefont {C.}~\bibnamefont {Desceliers}},\ }\href {https://doi.org/10.1007/s11831-021-09581-y} {\bibinfo {title} {{Review and Recent Developments on the Perfectly Matched Layer (PML) Method for the Numerical Modeling and Simulation of Elastic Wave Propagation in Unbounded Domains}}} (\bibinfo {year} {2022})\BibitemShut {NoStop}%
\bibitem [{\citenamefont {Skelton}\ \emph {et~al.}(2007)\citenamefont {Skelton}, \citenamefont {Adams},\ and\ \citenamefont {Craster}}]{Skelton2007GuidedLayers}%
  \BibitemOpen
  \bibfield  {author} {\bibinfo {author} {\bibfnamefont {E.~A.}\ \bibnamefont {Skelton}}, \bibinfo {author} {\bibfnamefont {S.~D.}\ \bibnamefont {Adams}},\ and\ \bibinfo {author} {\bibfnamefont {R.~V.}\ \bibnamefont {Craster}},\ }\bibfield  {title} {\bibinfo {title} {{Guided elastic waves and perfectly matched layers}},\ }\href {https://doi.org/10.1016/j.wavemoti.2007.03.001} {\bibfield  {journal} {\bibinfo  {journal} {Wave Motion}\ }\textbf {\bibinfo {volume} {44}},\ \bibinfo {pages} {573} (\bibinfo {year} {2007})}\BibitemShut {NoStop}%
\bibitem [{\citenamefont {Pelat}\ \emph {et~al.}(2011)\citenamefont {Pelat}, \citenamefont {Felix},\ and\ \citenamefont {Pagneux}}]{Pelat2011AWaveguides}%
  \BibitemOpen
  \bibfield  {author} {\bibinfo {author} {\bibfnamefont {A.}~\bibnamefont {Pelat}}, \bibinfo {author} {\bibfnamefont {S.}~\bibnamefont {Felix}},\ and\ \bibinfo {author} {\bibfnamefont {V.}~\bibnamefont {Pagneux}},\ }\bibfield  {title} {\bibinfo {title} {{A coupled modal-finite element method for the wave propagation modeling in irregular open waveguides}},\ }\href {https://doi.org/10.1121/1.3531928} {\bibfield  {journal} {\bibinfo  {journal} {J. Acoust. Soc. Am.}\ }\textbf {\bibinfo {volume} {129}},\ \bibinfo {pages} {1240} (\bibinfo {year} {2011})}\BibitemShut {NoStop}%
\bibitem [{\citenamefont {Bonnet-Ben~Dhia}\ \emph {et~al.}(2010)\citenamefont {Bonnet-Ben~Dhia}, \citenamefont {Goursaud}, \citenamefont {Hazard},\ and\ \citenamefont {Prieto}}]{Bonnet-BenDhia2010AWaveguides}%
  \BibitemOpen
  \bibfield  {author} {\bibinfo {author} {\bibfnamefont {A.~S.}\ \bibnamefont {Bonnet-Ben~Dhia}}, \bibinfo {author} {\bibfnamefont {B.}~\bibnamefont {Goursaud}}, \bibinfo {author} {\bibfnamefont {C.}~\bibnamefont {Hazard}},\ and\ \bibinfo {author} {\bibfnamefont {A.}~\bibnamefont {Prieto}},\ }\bibfield  {title} {\bibinfo {title} {{A multimodal method for non-uniform open waveguides}},\ }in\ \href {https://doi.org/10.1016/j.phpro.2010.01.065} {\emph {\bibinfo {booktitle} {Physics Procedia}}},\ Vol.~\bibinfo {volume} {3}\ (\bibinfo  {publisher} {Elsevier},\ \bibinfo {year} {2010})\ pp.\ \bibinfo {pages} {497--503}\BibitemShut {NoStop}%
\bibitem [{\citenamefont {Agha}\ \emph {et~al.}(2008)\citenamefont {Agha}, \citenamefont {Zolla}, \citenamefont {Nicolet},\ and\ \citenamefont {Guenneau}}]{Agha2008OnFibres}%
  \BibitemOpen
  \bibfield  {author} {\bibinfo {author} {\bibfnamefont {Y.~O.}\ \bibnamefont {Agha}}, \bibinfo {author} {\bibfnamefont {F.}~\bibnamefont {Zolla}}, \bibinfo {author} {\bibfnamefont {A.}~\bibnamefont {Nicolet}},\ and\ \bibinfo {author} {\bibfnamefont {S.}~\bibnamefont {Guenneau}},\ }\bibfield  {title} {\bibinfo {title} {{On the use of PML for the computation of leaky modes: An application to microstructured optical fibres}},\ }\href {https://doi.org/10.1108/03321640810836672} {\bibfield  {journal} {\bibinfo  {journal} {Int. J. Comput. Math. Electr.}\ }\textbf {\bibinfo {volume} {27}},\ \bibinfo {pages} {95} (\bibinfo {year} {2008})}\BibitemShut {NoStop}%
\end{thebibliography}%

\end{document}